\documentclass[12pt,a4paper,oneside]{article}
\usepackage{a4wide}
\usepackage[T1]{fontenc}
\usepackage[utf8]{inputenc}
\usepackage{amsmath}
\usepackage{times}
\usepackage{graphicx}
\usepackage[final]{pdfpages}
\usepackage{longtable}
\usepackage{float}
\usepackage{amssymb,bm}
\usepackage{url}
\usepackage{xcolor}
\usepackage{mathtools}
\usepackage{tcolorbox}
\usepackage{multirow}
\definecolor{pigpink}{HTML}{FDD7E4}
\definecolor{lcyan}{HTML}{E0FFFF}
\definecolor{mint}{HTML}{98FF98}

\usepackage{marginnote}
\usepackage{soul}

\usepackage{amsthm}
\newtheorem{theorem}{Theorem}[section]

\newtheorem{definition}[theorem]{Definition}
\newtheorem{remark}[theorem]{Remark}

\newcommand{\Id}{\mathrm{Id}}
\newcommand{\Q}{\mathcal{Q}}

\newcommand{\bs}{\boldsymbol{\sigma}}
\newcommand{\Lo}{\mathcal{L}}

\usepackage{yfonts}

\title{
Efficient Numerical Algorithms for Phase-Amplitude Reduction on the Slow Attracting Manifold of Limit cycles}

\author{David Reyner-Parra$^{1}$, Alberto P\'erez-Cervera$^{1}$ \& Gemma Huguet$^{1,2,3}$ \\
\parbox{12.5cm}{
  \small
  \begin{itemize}
  \item[$^1$]
    Departament de Matem\`atiques, Universitat Polit\`ecnica de Catalunya, Barcelona, Spain 
  \item[$^2$]
    Institut de Matem\`atiques de la UPC - Barcelona Tech (IMTech), Barcelona, Spain 
  \item[$^ 3$]
  Centre de Recerca Matem\`atica, Barcelona, Spain
   \end{itemize}
 }}
 \date{}

\begin{document}
\maketitle

\textbf{Corresponding author:} Gemma Huguet, \texttt{gemma.huguet@upc.edu} \\


\begin{abstract}
The phase-amplitude framework extends the classical phase reduction method by incorporating amplitude coordinates (or isostables) to describe transient dynamics transverse to the limit cycle in a simplified form. While the full set of amplitude coordinates provides an exact description of oscillatory dynamics, it maintains the system's original dimensionality, limiting the advantages of simplification. A more effective approach reduces the dynamics to the slow attracting invariant submanifold associated with the slowest contracting direction, achieving a balance between simplification and accuracy. In this work, we present an efficient numerical method to compute the parameterization of the attracting slow submanifold of hyperbolic limit cycles and the simplified dynamics in its induced coordinates. Additionally, we compute the infinitesimal Phase and Amplitude Response Functions (iPRF and iARF, respectively) restricted to this manifold, which characterize the effects of perturbations on phase and amplitude. These results are obtained by solving an invariance equation for the slow manifold and adjoint equations for the iPRF and iARF. 
To solve these functional equations efficiently,  we employ the Floquet normal form to solve the invariance equation and propose a novel coordinate transformation to simplify the adjoint equations. The solutions are expressed as Fourier-Taylor expansions with arbitrarily high accuracy. Our method accommodates both real and complex Floquet exponents. Finally, we discuss the numerical implementation of the method and present results from its application to a representative example.
\end{abstract}

\section{Introduction}

Oscillations are ubiquitous across physical, chemical, and biological systems, appearing in diverse processes \cite{kuramoto2003chemical, winfree2001geometry, strogatzbook, PIK01}. From a dynamical systems perspective, these oscillations correspond to attracting limit cycles in phase space. The dynamics on the limit cycle is effectively captured by a phase variable. A widely used approach for analyzing perturbed oscillators and systems of coupled oscillators is the phase reduction \cite{hoppensteadt2012, ErmentroutTerman2010, ErmentroutKopell91, winfree1967biological}. This method assumes that the dynamics of the perturbed system remains near the unperturbed limit cycle, with perturbations affecting only the dynamics along the cycle. Phase reduction offers significant simplification by reducing the system of ODEs to a single equation, while preserving the most essential dynamical properties. However, its applicability is constrained by the strong assumption that the dynamics are restricted to a neighborhood of the unperturbed limit cycle, which may not always hold. 

To overcome these limitations, recent research has focused on incorporating additional variables to describe dynamics along directions transverse to the limit cycle. Various approaches to define this transverse coordinate system have been proposed by different groups \cite{wedgwood2013phase, shirasaka2017phase, bonnin17, mauroy2018global, moehliswilsonpre2016, wilson2018greater, Wilson2020}. Of particular interest for this paper are the works in \cite{guillamon2009computational, huguet2013computation, castejon2013phase, perezrole, PerezCervera2020}, which build upon the parameterization method introduced in \cite{cabre2003parameterization, cabre2003parameterization2, cabre2005parameterization} (for a review, see \cite{haro2016}). This method parameterizes the stable invariant manifold of the hyperbolic limit cycle using phase-amplitude variables, which are designed to follow a prescribed simple dynamics. This setup is formalized through a functional equation, known as the invariance equation, which links the parameterization to the intrinsic dynamics. For a hyperbolic attracting limit cycle, the stable manifold coincides with its basin of attraction. Thus, the parameterization method provides a coordinate transformation to phase-amplitude variables, for which the dynamics has a simpler expression. Furthermore, this parameterization allows to define the infinitesimal Phase and Amplitude Response Functions (iPRF and iARF), which characterize how infinitesimal perturbations affect the phase and amplitude variables, respectively, along the trajectories \cite{guillamon2009computational, castejon2013phase, Wilson2020}. These functions allow for the study of perturbed dynamics within the phase-amplitude framework.

Using the full set of amplitude coordinates is generally not advantageous, as it retains the original system's dimensionality despite simplifying the dynamics. A more practical reduction strategy, less restrictive than phase reduction but still effective, is to focus on the dynamics restricted onto the slow attracting invariant submanifold (assuming it exists), which corresponds to the manifold associated to the smallest (in modulus) Floquet exponent \cite{cabre2003parameterization}. The advantage of this approach is that these reduced systems can be incorporated into larger networks of oscillators, offering more accurate results and capturing richer dynamics than the simple phase description \cite{Nicks2024}, while minimally impacting the system's dimensionality growth. 

To perform this reduction it is necessary to know the parameterization of the slow submanifold, as well as the functions iPRF and iARF on this manifold to describe the effects of perturbations. Certainly, computing the full parameterization of the stable manifold automatically yields the slow submanifold and the dynamics on it, as it can be obtained by setting the remaining amplitude coordinates to zero within the parameterization. However, computing the full parameterization can be computationally expensive, especially for high-dimensional systems. 

In this paper, we introduce an efficient numerical method to compute both the slow submanifold and the iPRF and the iARF functions on it, without computing the full parameterization. Our approach uses the invariance equation for the slow submanifold and the adjoint equations for the iPRF and the iARFs as introduced in \cite{castejon2013phase} (and further explored in \cite{Wilson2020}) and presents a strategy for solving them efficiently. Indeed, we impose a formal power series solution for these equations and by matching terms with the same power we obtain a sequence of linear ordinary differential equations with periodic coefficients, which need to be solved recursively. Using the Floquet normal form \cite{castelli2015parameterization, PerezCervera2020, huguet2013computation}, we can transform each differential equation in the recursive scheme into a non-homogeneous linear one with constant coefficients. This approach simplifies the problem into algebraic equations that are diagonal in Fourier space. As a novel contribution, we extend this methodology to efficiently solve the adjoint equations. Specifically, we choose an appropriate coordinate system consisting of the infinitesimal Phase Response Curve (iPRC) and the infinitesimal Amplitude Response Curve (iARC), which are the iPRF and iARF functions evaluated on the limit cycle. As in the invariance equation, this transformation leads to simple equations for the Fourier-Taylor coefficients of the solutions of the adjoint equations. 

As in previous works \cite{PerezCervera2020, huguet2013computation}, we use automatic differentiation algorithms \cite{griewank2008evaluating, haro2016} to compute compositions of power series with elementary functions, thus avoiding the computation of high-order derivatives of the vector field.
Additionally, our method accommodates both real and complex Floquet exponents, expanding upon the cases considered in \cite{PerezCervera2020} and \cite{Wilson2020}.
 
We illustrate the techniques by applying them to a 6-dimensional model from \cite{dumont2019macroscopic} that describes the exact mean-field dynamics of a neural network of excitatory and inhibitory cells (E-I network). The system represents a comprehensive example of most of the remarkable key points of the methodology, involving complex Floquet exponents. Furthermore, the model and variants of it has been widely used in the literature to study several aspects related to neural oscillations \cite{Segneri2020,Clusella2023,ReynerHuguet22, Orieux2024, Taher2020}.

The structure of the paper is as follows. In Section \ref{sec:section2} we introduce the theoretical and computational methodology to obtain the parameterization of the invariant manifold as well as the iPRF and iARFs. In Section \ref{sec:sec-3} we present the numerical strategy to solve the adjoint equations efficiently, and obtain a parameterization of the slow submanifold and the dynamics restricted to it. In Section \ref{sec:numerics} we present numerical results for a specific system widely studied in computational neuroscience, along with a discussion of implementation details for the algorithms. Finally, we conclude with the discussion in Section \ref{sec:discussion}, where we review our findings and their relation to other results in the field. The Appendix is dedicated to demonstrating the existence and uniqueness of solutions for the adjoint equations involved in the computation of the slow reduction dynamics.

\section{Background and statement of the problem}\label{sec:section2}

In this section, we review the basic principles of the parameterization method used to describe oscillators in phase-amplitude coordinates and its response to perturbations. 

\subsection{Phase-Amplitude variables and the parameterization method} 
Let us consider an autonomous system of ODEs 
\begin{align}\label{eq:mathDef_1}
    \dot{x} = X(x), \quad x \in \mathbb{R}^{d}, \quad d \geq 2 ,
\end{align}
whose flow is denoted by $\phi_t(x)$. Here, and throughout the manuscript, $\dot{x}$ denotes the derivative with respect to time. Assume that $X: \mathbb{R}^d \rightarrow \mathbb{R}^d$ is an analytic vector field and that system \eqref{eq:mathDef_1} has a $T$-periodic hyperbolic attracting limit cycle $\Gamma$, parameterized by the phase variable $\theta = t/T$ through the function
\begin{equation}
    \begin{aligned}\label{eq:mathDef_2}
        \gamma:\mathbb{T}:= \mathbb{R}/\mathbb{Z} &\to \mathbb{R}^{d}\\
        \theta &\mapsto \gamma(\theta).
    \end{aligned}
\end{equation}
The function $\gamma$ has period 1, that is, $\gamma(\theta) = \gamma(\theta +1)$ and $x(t) = \gamma(\frac{t}{T})$ satisfies \eqref{eq:mathDef_1}. Thus, the dynamics of \eqref{eq:mathDef_1} on $\Gamma$ can be reduced to a single equation for the phase variable
\begin{equation}\label{eq:phase-reduction}
    \dot{\theta} = \frac{1}{T}, \quad \quad \theta \in \mathbb{T}.
\end{equation}

Since we are considering a hyperbolic attracting periodic orbit, all its characteristic exponents have negative real part, except the trivial one which is 0 (or, equivalently, the Floquet multipliers are inside the unit circle except the trivial one which is 1). The characteristic exponents of $\Gamma$ can be obtained by solving the variational equations of system \eqref{eq:mathDef_1} along the solution $\gamma(t/T)$, that is,
\begin{equation}\label{eq:mjVarEqs}
    \dot{\Phi} = DX \big(\gamma(t/T) \big) \Phi, \quad \quad \quad \text{with} \quad \Phi(0) = \Id_{d \times d} \, .
\end{equation}
The solution $\Phi(t)$ of the variational equations \eqref{eq:mjVarEqs} evaluated at $t = T$, i.e. $\Phi_T:=\Phi(T)$, is known as the monodromy matrix. Then, the eigenvalues of $\Phi_T$, namely $\mu_j$, $j = 0,\ldots, d-1$, are the \emph{Floquet or characteristic multipliers} of the limit cycle $\Gamma$ and the values $\lambda_j \in \mathbb{C}$ such that $ \mu_j= e^{\lambda_j T}$ are the \emph{Floquet or characteristic exponents}. 
We label the Floquet exponents so that $\Re(\lambda_{j+1})\leq \Re(\lambda_j)$, for $j=0,\ldots,d-1$. The index $j = 0$ will be assigned to the trivial multiplier $\mu_0 = 1$, so $\lambda_0 = 0$. 

\begin{remark}
The Floquet multipliers are uniquely determined. However, the Floquet exponents are undetermined except for the sum of multiples of $2 \pi i/T$. The unique real numbers $\ell_j$ such that $|\mu_j|=e^{\ell_{j }T}$ are called \emph{Lyapunov exponents}.
\end{remark}

The parameterization $\gamma$ in \eqref{eq:mathDef_2} of the limit cycle $\Gamma$ together with the restricted dynamics provided by equation \eqref{eq:phase-reduction} are the centerpiece of the so-called phase reduction method for oscillators. The parameterization method \cite{cabre2003parameterization, cabre2003parameterization2, cabre2005parameterization} shares a similar objective with the phase reduction: it seeks a parameterization of an $n$-dimensional invariant manifold of a dynamical system and expresses the dynamics on it in the coordinates induced by this parameterization. These coordinates are chosen so that the dynamics are as simple as possible, often rendering them linear. When the manifold considered is the stable manifold of an attracting limit cycle, it coincides with the basin of attraction $\Omega \subset \mathbb{R}^ d$ of the limit cycle $\Gamma$ (see \cite{guckenheimer1975}). Thus, this is equivalent to finding a change of variables $K$ that conjugates the vector field $X$ in $\Omega \subset \mathbb{R}^ d$ to a vector field $\mathcal{X}$ with a simpler expression of the  dynamics (\cite{guillamon2009computational, huguet2013computation, PerezCervera2020}). In this paper, we will consider the parameterization of the slow submanifold of the stable manifold, associated to the non-trivial Floquet exponent of smallest modulus (as we will see in Section \ref{sec:sec-3}).

The parameterization method (see the book \cite{haro2016} for a review) states a functional equation for the parameterization that characterizes the invariance of the manifold and the dynamics on it in the coordinates defined by the parameterization. Thus, we look for a local analytic diffeomorphism
\begin{equation}\label{eq:kThetaSigma}
\begin{aligned}
K : \mathbb{T} \times (\mathcal{B} \subset \mathbb{C}^{d-1}) &\rightarrow \mathbb{R}^d \\
(\theta, \bs) &\rightarrow K(\theta, \bs) \, , 
\end{aligned}
\end{equation}
that satisfies the invariance equation
\begin{equation}\label{eq:invariantgen}
    DK(\theta,\bs) \cdot \mathcal{X}_{\Lambda} = X \circ K(\theta,\bs) \, ,
\end{equation}
where 
\[
\mathcal{X}_{\Lambda}:=
\begin{pmatrix}
    1/T \\ \Lambda \, \bs 
\end{pmatrix}
\quad \textrm{with} \quad
\Lambda :=
\begin{pmatrix}
    \lambda_1 & & \\
    & \ddots & \\
    & & \lambda_{d-1}  
\end{pmatrix} \, .
\]
Here, $T$ is the period of the limit cycle $\Gamma$ and $\lambda_1, \ldots, \lambda_{d-1} \in \mathbb{C}$ its non-trivial characteristic exponents. 
We can also write \eqref{eq:invariantgen} as
\begin{equation}\label{eq:mjInvEq}
    \frac{1}{T}\frac{\partial}{\partial\theta}K(\theta, \bs) + \sum_{j=1}^{d-1} \lambda_j \sigma_j \frac{\partial}{\partial \sigma_j}K(\theta, \bs) = X(K(\theta, \bs)) \, .
\end{equation}

The coordinates of the parameterization $K$ are the phase $\theta \in \mathbb{T}$ introduced in \eqref{eq:mathDef_2} and the amplitude coordinates $\bs=(\sigma_1,\ldots,\sigma_{d-1}) \in \mathbb{C}^{d-1}$, corresponding to transverse directions to the limit cycle. Notice that we do not assume that the Floquet exponents are real, and for this reason the amplitude variable $\bs \in \mathbb{C}^{d-1}$. The dynamics of the vector field $X$ in \eqref{eq:mathDef_1} expressed in these new variables $(\theta, \bs) \in \mathbb{T} \times (\mathcal{B} \subset \mathbb{C}^{d-1})$ is given by 
\begin{equation}\label{eq:phase-amplitude}
    \dot{\theta} = \frac{1}{T}, \quad \quad \quad \dot{\bs} = \Lambda \bs \, ,
\end{equation} 
that is, the variable $\theta$ rotates at a constant speed $1/T$, while each variable $\sigma_j \in \mathbb{C}$ contracts at a rate $\Re(\lambda_j)$. Thus, the evolution of the flow $\phi_t$ in the coordinates $(\theta,\bs)$  becomes
\begin{equation}\label{eq:aboveEq}
    \phi_t(K(\theta, \bs)) = K\Big(\theta + \frac{t}{T}, e^{\Lambda t}\bs\Big) \, .
\end{equation}

It is known (see \cite{cabre2005parameterization}) that provided the characteristic exponents satisfy certain non-resonance conditions, one can indeed solve the functional equation \eqref{eq:mjInvEq} and find a map $K(\theta, \bs)$, at least formally. A formal series solution of the invariance equation~\eqref{eq:mjInvEq} following \cite{castelli2015parameterization, PerezCervera2020} is discussed in Section \ref{sec:computation_K1}.

We introduce here the definition of resonances in our setting.

\begin{definition}\label{def:res}
We say that the Floquet exponents $(\lambda_1,\ldots,\lambda_{d-1})$ are resonant if there exists a ($d-1$)-dimensional multi-index $a=(a_1,\ldots,a_{d-1}) \in \mathbb{N}^{d-1}$ (non-negative integers) and an index $k \in \{1,\ldots,d-1\}$ such that
\[\sum_{i=1}^{d-1} a_i \lambda_i - \lambda_k=0, \quad \textrm{with} \quad |a| :=a_1+\cdots+a_{d-1} \geq 2.\]
We say that $|a|$ is the order of the resonance.
\end{definition}

\begin{remark}\label{rem:nu}
Notice that equation \eqref{eq:mjInvEq} does not have a unique solution. Indeed, if $K(\theta, \bs)$ is a solution then $K(\theta + \theta_0, \mathbf{b} \, \bs)$ is also a solution, for any $\theta_0 \in \mathbb{T}$ and $\mathbf{b} \in \mathbb{R}^{d-1}$. The meaning of $\theta_0$ is the choice of the origin of time, and $b_j$ corresponds to the choice of units in $\sigma_j$.
\end{remark}

The map $K$ in \eqref{eq:kThetaSigma} allows to define a scalar function $\Theta$ that assigns the asymptotic phase to any point $x$ in the basin of attraction $\Omega \subset \mathbb{R}^d$ of the limit cycle $\Gamma$. Indeed, 
\begin{equation}
    \begin{aligned}
        \Theta:\Omega \subset \mathbb{R}^{d} &\to \mathbb{T},\\
        x &\mapsto \Theta(x) = \theta \quad \quad \text{where} \enskip x = K(\theta, \bs),\enskip \textrm{for some} \enskip \bs \in \mathbb{R}^{d-1} \, .
    \end{aligned}\label{eq:mathDef_3}
\end{equation}
Notice that $\Theta(\phi_t(x)) = \Theta(x) + \frac{t}{T}$. The level curves of $\Theta$ correspond to the isochrons $\mathcal{I}_{\theta}$ (\cite{winfree1967biological, guckenheimer1975}), that is,
\begin{equation}
    \mathcal{I}_{\theta} = \{x \in \Omega \subset \mathbb{R}^d \enspace | \enspace \Theta(x) = \theta \} \, .
\end{equation}

Analogously, the map $K$ in \eqref{eq:kThetaSigma} also allows us to define the scalar complex-valued functions $\Sigma_j$, for $j=1,\ldots, d-1$ that assign the amplitude variable $\sigma_j$ to any point $x \in$ $\Omega$:
\begin{equation}\label{eq:sigma3Dcase}
    \begin{aligned}
        \Sigma_j : \Omega \subset \mathbb{R}^{d} &\to \mathbb{C},\\
        x &\mapsto \Sigma_j(x) = \sigma_j, \quad \text{where} \enskip x = K(\theta, \bs), \enskip \textrm{for some} \enskip \theta \in \mathbb{T}.
    \end{aligned}
\end{equation}
Notice that $\Sigma_j(\phi_t(x)) = \Sigma_j(x)e^{\lambda_j t}$, for $j=1,\ldots,d-1$ . The sets with the same asymptotic convergence to the limit cycle are called \textit{isostables} (see \cite{mauroy2018global, mauroy2013isostables}) or \emph{A-curves} for the 2-dimensional case (see \cite{castejon2013phase}), and correspond to the sets of points 
\begin{equation}\label{eq:aCurvesDef}
    \mathcal{A}_{s}^j = \{x \in \Omega \subset \mathbb{R}^d \enspace | \enspace |\Sigma_j(x)| = s \},
\end{equation}
where $s \in \mathbb{R}$ and $|\cdot|$ denotes the modulus. We will denote by $\Sigma$ the vector-valued function $\Sigma(x):=(\Sigma_1(x),\ldots,\Sigma_{d-1}(x))$.

Finally, we recall that the vector-valued
function $(\Theta,\Sigma)$ is the inverse of $K$ since
\begin{equation}\label{eq:Kinv}
    K\big(\Theta(x), \Sigma_1(x), \ldots, \Sigma_{d-1}(x)\big) = x \, , \quad \text{for} \enskip x = (x_1, \ldots, x_d) \in \Omega \subset \mathbb{R}^d.
\end{equation}

\subsection{Normal bundle to the periodic orbit}\label{sec:computation_K1}
In this section we discuss how to obtain a formal series solution of the invariance equation~\eqref{eq:mjInvEq}. As we will see in Section~\ref{parameterization_slow_manifold}, for the main algorithm presented in this work, we only need the normal bundle of the periodic orbit. For this reason, we only discuss the solutions of the resulting equations up to order 1.

We assume a formal series solution for the parameterization $K$ in equation \eqref{eq:mjInvEq} of the form
\begin{equation}\label{eq:formal}
    K(\theta,\bs)=\sum_{|\alpha|=0}^{\infty} K_{\alpha}(\theta) \bs^{\alpha},
\end{equation}
where $\alpha \in \mathbb{N}^{d-1}$ is a $(d-1)$-dimensional multi-index, $|\alpha| :=\alpha_1+\cdots+\alpha_{d-1} $ and $\bs^{\alpha}:=\sigma_1^{\alpha_1} \ldots \sigma_{d-1}^{\alpha_{d-1}}$. 
We substitute \eqref{eq:formal} into equation \eqref{eq:mjInvEq} and collect the terms with the same power of $\bs$, thus obtaining the following set of differential equations:
\begin{itemize}
    \item Order $|\alpha|=0$. The only multi-index is $\mathbf{0}=(0,\ldots,0)$, the term $\bs^0$ gives the following equation for the 1-periodic function $K_{\mathbf{0}}(\theta)$:
    \begin{equation}\label{eq:eqMizero}
        \frac{1}{T}\frac{d}{d\theta}K_{\mathbf{0}}(\theta) = X(K_{\mathbf{0}}(\theta)),
    \end{equation}
    The solution is given by the periodic orbit, i.e. $K_{\mathbf{0}}(\theta)=\gamma(\theta)$.

    \item Order $|\alpha|=1$. The multi-indices of length 1 are $e_j=(0,\ldots,1,\ldots,0)$ with 1 in the $j$-th position. We obtain the following equation for the 1-periodic functions $K_{e_j}(\theta)$ and scalars $\lambda_j \in \mathbb{C}$:
    \begin{equation}\label{eq:eqMione}
        \frac{1}{T}\frac{d}{d\theta}K_{e_j}(\theta) + \lambda_i K_{e_j}(\theta)= DX(K_{\mathbf{0}}(\theta)) K_{e_j}(\theta).
    \end{equation}
    
     The solutions of equations \eqref{eq:eqMione} (see \cite{cabre2005parameterization, castelli2015parameterization, huguet2013computation}) are given by
    \begin{equation}\label{eq:K1sol_complex}
        K_{e_j}(\theta)= e^{-\lambda_j T \theta}\Phi(T \theta)w_j \, ,
    \end{equation} 
    where $(e^{\lambda_j T},w_j)$ is an eigenpair of the monodromy matrix $\Phi(T)$ of the variational equations \eqref{eq:mjVarEqs}. Notice that when $\lambda_j$ is complex so is the corresponding function $K_{e_j}(\theta)$ as well as its associated amplitude variable $\sigma_j$.

    \item Order $|\alpha| \geq 2$, the equations for the 1-periodic functions $K_{\alpha}(\theta)$ are of the form
    \begin{equation}\label{eq:eqMihorder}
    \frac{1}{T}\frac{d}{d\theta}K_{\alpha}(\theta) + (\alpha_1 \lambda_1 + \ldots + \alpha_{d-1} \lambda_{d-1}) K_{\alpha}(\theta)= DX(K_{\mathbf{0}}(\theta)) K_{\alpha}(\theta) + B_{\alpha}(\theta) \, ,
    \end{equation}
    where $B_{\alpha}$ is a polynomial that involves only lower order terms. It is possible to design strategies to solve these equations efficiently by means of the Floquet normal form (see \cite{PerezCervera2020, castelli2015parameterization}). However, we do not discuss them in this paper, since they are not necessary for the purposes of the paper.
\end{itemize}

\subsection{Real parameterization of the invariant manifold}
Up to this point, we have considered the case in which the matrix $\Lambda$ and the amplitude variables $\bs$ in the parameterization $K$ are complex. However, it is possible to consider real amplitude variables $\widetilde \bs \in \mathbb{R}^{d-1}$. In this case, we introduce the real parameterization 
\begin{equation}\label{eq:kThetaSigma_real}
    \begin{aligned}
        \widetilde K : \mathbb{T} \times (\widetilde{\mathcal{B}} \subset \mathbb{R}^{d-1}) &\rightarrow \mathbb{R}^d \\
        (\theta, \widetilde \bs) &\rightarrow \widetilde K(\theta, \widetilde \bs) \, , 
    \end{aligned}
\end{equation}
satisfying the following invariance equation
\begin{equation}\label{eq:inveq_real}
    D \widetilde K(\theta, \widetilde \bs) \cdot \mathcal{X}_{\Lambda_R} = X \circ \widetilde K(\theta, \widetilde \bs)
\end{equation}
where 
\begin{equation}\label{eq:Lambda_real}
    \mathcal{X}_{\Lambda_R}:=
    \begin{pmatrix}
        1/T \\ \Lambda_R \, \widetilde \bs
    \end{pmatrix}
    \quad \textrm{and} \quad 
    \Lambda_R:=
    \begin{pmatrix}
        R_1 & & \\
        & \ddots & \\
        & & R_{d_m}  
    \end{pmatrix} \, ,
\end{equation}
with $d_m=d-1-m$, where $m$ is the number of pairs of complex conjugates Floquet multipliers. The blocks $R_j$ are given by:
\begin{itemize}
    \item If $\lambda_j \in \mathbb{R}$ (real positive Floquet multiplier $\mu_j \in \mathbb{R}^+$), then $R_j=(\lambda_j)$ and $\widetilde \sigma_j=\sigma_j \in \mathbb{R}$.

    \item If $\lambda_{j}=\nu_j + i\pi/T=\dfrac{1}{T} (\ln \, |\mu_j|+ i \pi)$ (real negative Floquet multiplier $\mu_j<0$), then
    $R_j=(\nu_j)$ and $\widetilde \sigma_j=|\sigma_j| \in \mathbb{R}$. The parameterization satisfies
    \[\widetilde K (\theta, 0, \ldots, \widetilde \sigma_j,\ldots,0)= \widetilde K (\theta+1, 0, \ldots, -\widetilde \sigma_j,\ldots,0).\]

    \item If $\lambda_{j,j+1}=\alpha_j \pm i \beta_j$ are complex conjugates (complex conjugate Floquet multipliers $\mu_{j,j+1} \in \mathbb{C}$) then
    $R_{j}=\begin{pmatrix}
    \alpha_j & \beta_j \\
    -\beta_j & \alpha_j \\  
    \end{pmatrix}$ and $R_{j+1}=\varnothing$. The pair of complex conjugates amplitude variables $\sigma_{j}=\overline{\sigma}_{j+1} \in \mathbb{C}$ is related to the two real amplitude variables $\widetilde \sigma_j, \widetilde \sigma_{j+1} \in \mathbb{R}$ by 
    \begin{equation}\label{eq:sigma_cr}
        \widetilde \sigma_j=\sigma_j + \overline \sigma_j=2 \Re\{\sigma_j\} \quad \textrm{and} \quad \widetilde \sigma_{j+1}=i(\sigma_j - \overline \sigma_j)=-2 \Im\{ \sigma_j\} \, .
    \end{equation}
\end{itemize}

Therefore, we can define the real amplitude functions $\widetilde \Sigma_{j}: \Omega \subset \mathbb{R}^d \rightarrow \mathbb{R}$, satisfying $\widetilde \Sigma_{j}(\widetilde K(\theta,\widetilde \bs))=\widetilde \sigma_{j}$ which are related with the complex amplitude functions $\Sigma_{j}$ by
\begin{itemize}
    \item If $\lambda_j \in \mathbb{R}$, then $\widetilde\Sigma_j(x)=\Sigma_j (x)$.
    
    \item If $\lambda_j =\nu + i\pi/T \in \mathbb{C}$, then 
    \begin{equation}\label{eq:icomplexa}
        \widetilde \Sigma_j(x)= |\Sigma_j (x)|.
    \end{equation}
    
    \item If $\lambda_{j,j+1}=\alpha \pm i \beta \in \mathbb{C}$, then 
    \begin{equation}\label{eq:real-and-complex-amplitudes}
        \begin{split}
            \widetilde \Sigma_j(x) &= \Sigma_j (x) + \Sigma_{j+1} (x) = 2 \Re \{\Sigma_j (x)\}, \\
            \widetilde \Sigma_{j+1}(x) &= i(\Sigma_j (x)- \Sigma_{j+1} (x)) = -2 \Im \{\Sigma_j (x)\}.
        \end{split}
    \end{equation}
\end{itemize}

From expression \eqref{eq:formal}, it is possible to obtain a formal series expansion for the real parameterization $\widetilde K$ satisfying equation \eqref{eq:inveq_real}. Indeed, if we denote
\begin{equation}\label{eq:Kformal_real}
    \widetilde K (\theta, \widetilde \bs)= \sum_{|\alpha|=0}^{\infty} \widetilde K_{\alpha}(\theta) \widetilde \bs^{\alpha},
\end{equation}
then using expressions in \eqref{eq:sigma_cr} and collecting terms with the same power in $\tilde \bs$, we can obtain the terms $\widetilde K_{\alpha}$ from $K_{\alpha}$. 

However, it is possible to obtain the functions $\widetilde K_{\alpha} (\theta)$ directly using the same strategy as described before, but in this case looking for a solution in the form \eqref{eq:Kformal_real} of the invariance equation \eqref{eq:inveq_real}. Next, we discuss how to solve the resulting equations only up to order 1 since, as mentioned before, higher order terms are not needed for the purposes of this paper.

Of course, the equation for $\widetilde K_{\mathbf{0}}$ is the same as for $K_{\mathbf{0}}$ and $\widetilde K_{\mathbf{0}}(\theta) = K_{\mathbf{0}} (\theta)=\gamma(\theta)$. The equations for the terms $\widetilde K_{e_j}$ of order 1 are, depending on the Floquet exponent, 
\begin{itemize}
    \item If $\lambda_j \in \mathbb{R}$, then $\widetilde K_{e_j}(\theta)$ is a 1-periodic real function satisfying
    \[\frac{1}{T}\frac{d}{d\theta} \widetilde K_{e_j}(\theta) + \lambda_j \widetilde K_{e_j}(\theta)= DX(K_{\mathbf{0}}(\theta)) \widetilde K_{e_j}(\theta).\]
    In this case, $\widetilde K_{e_j}(\theta)=K_{e_j}(\theta)$.

    \item If $\lambda_{j}=\nu_j + i \dfrac{\pi}{T}$, $\nu=\ln |\mu_j| \in \mathbb{R}$, then $\widetilde K_{e_j}(\theta)$ is a 2-periodic real function satisfying
    \[\frac{1}{T} \frac{d}{d \theta} \widetilde K_{e_i}(\theta) + \nu_j \widetilde K_{e_i}(\theta) = DX (K_{\mathbf{0}}(\theta)) \widetilde K_{e_i} (\theta). \]
    The solution is given by 
    \begin{equation}\label{eq:Kim}
    \widetilde K_{e_j}(\theta)=e^{-\nu_j T \theta} \Phi(T \theta) w_j,
    \end{equation}
    where recall that the eigenvector $w_j \in \mathbb{R}^d$. From the expression \eqref{eq:Kim}, it is clear that $\widetilde K_{e_j}(\theta+1)=- \widetilde K_{e_j}(\theta)$, so $\widetilde K_{e_j}$ is $2$-periodic. Thus, we restrict to $\widetilde \Sigma_j(x)=\widetilde \sigma_j>0$ and $\widetilde K_{e_j}$ defined for $\theta \in [0,2)$. In this case $K_{e_j} (\theta)=e^{-i \pi \theta} \widetilde K_{e_j}(\theta)$.

    \item If $\lambda_{j,j+1}=\alpha_j \pm i\beta_j$ are complex conjugates, with $\alpha_j,\beta_j \in \mathbb{R}$ then
    $\widetilde K_{e_j}$ and $\widetilde K_{e_{j+1}}$ are 1-periodic real functions satisfying the system of equations
    \begin{equation}
        \begin{array}{rcl}
            \dfrac{1}{T} \dfrac{d}{d \theta} \widetilde K_{e_j}(\theta) +  \alpha_j \widetilde K_{e_j}(\theta) - \beta_j \widetilde K_{e_{j+1}}(\theta)  &=& DX (K_{\mathbf{0}}(\theta)) \widetilde K_{e_j} (\theta), \\[0.75em]
            \dfrac{1}{T} \dfrac{d}{d \theta} \widetilde K_{e_{j+1}}(\theta) + \beta_j \widetilde K_{e_j}(\theta) + \alpha_j \widetilde K_{e_{j+1}}(\theta)  &= &DX (K_{\mathbf{0}}(\theta)) \widetilde K_{e_{j+1}} (\theta). 
        \end{array}
    \end{equation}
    So, let $w_{j}(\theta)=\Phi(\theta T) w_{j} = r_j(\theta) + i s_j(\theta)$, where the eigenvector $w_j$ is complex, then
    \begin{equation}\label{eq:Kpairreal}
    \begin{array}{rl}
       \widetilde K_{e_j}(\theta)&= e^{-\alpha_j \theta T} (r_j(\theta) \cos (\beta_j \, \theta T) + s_j(\theta) \sin (\beta_j \, \theta T) ) \, , \\
       \widetilde K_{e_{j+1}}(\theta)&= e^{-\alpha_j \theta T} (s_j(\theta) \cos (\beta_j \, \theta T) - r_j(\theta) \sin (\beta_j \, \theta T) ) \, . \\
    \end{array}
    \end{equation}
    Notice that expressions in \eqref{eq:Kpairreal} are just the real and imaginary parts of the complex solution \eqref{eq:K1sol_complex} of equation \eqref{eq:eqMione} with $\lambda_j=\overline{\lambda}_{j+1}$. Therefore, the relationship with $K_{e_j}$ is given by
    \[\widetilde K_{e_j}(\theta)=\Re \{K_{e_j}(\theta) \}, \qquad \widetilde K_{e_{j+1}}(\theta)=\Im \{ K_{e_j}(\theta) \}.\]
\end{itemize}

\subsection{Perturbed Oscillatory Dynamics in Phase-Amplitude coordinates}\label{sec:phaseAmpFun}
Consider the perturbed system of ODEs
\begin{equation}\label{eq:syst_pert}
    \dot{x}=X(x) + G(x,t) \, ,
\end{equation}
where $G:\Omega \subset \mathbb{R}^d \times \mathbb{R} \rightarrow \mathbb{R}^d$ is an analytic function. Then, by means of the change of coordinates $K$ defined in \eqref{eq:kThetaSigma}, the dynamics of the perturbed system
\eqref{eq:syst_pert} in the phase-amplitude coordinates is given by
\begin{equation}
    \begin{split}
        \dot{\theta} &=\dfrac{1}{T}+\nabla \Theta(K(\theta,\bs)) \cdot G(K(\theta,\bs),t) \, , \\[0.5em]
        \dot{\sigma}_j &=\lambda_j \sigma_j +\nabla \Sigma_j(K(\theta,\bs))\cdot G(K(\theta,\bs),t) \, , \quad \mbox{ for } j=1,\ldots,d-1 \, , \\
    \end{split}
\end{equation}
where $\cdot$ denotes the dot product in $\mathbb{R}^d$. The functions $\nabla\Theta(x)$ and $\nabla\Sigma_j(x)$ are the so-called infinitesimal Phase Response Function (iPRF) and the infinitesimal Amplitude Response Functions (iARFs), respectively (see \cite{Wilson2020, castejon2013phase, PerezCervera2020}). They quantify the shift in phase and amplitude when a perturbation is applied at a point $x$. Indeed, when these functions are restricted to the limit cycle, they become the classical infinitesimal Phase Response Curve (iPRC) and infinitesimal Amplitude or Isostable Response Curves (iARC) \cite{ErmentroutTerman2010}, that is,
\[\text{iPRC}(\theta)=\nabla \Theta (\gamma(\theta)) \quad \mbox{  and } \quad
\text{iARC}_j(\theta)=\nabla \Sigma_j (\gamma(\theta)) \, .\]

The iPRF $\nabla\Theta(x)$ and the iARFs $\nabla\Sigma_j(x)$ can be computed by means of the parameterization $K(\theta, \bs)$. Indeed, taking derivatives on both sides at the expression \eqref{eq:Kinv} we have 
\begin{equation}\label{eq:usefulSpace} 
    \begin{bmatrix} 
        DK(\theta, \bs) 
    \end{bmatrix}
    \begin{bmatrix} 
        \big\uparrow & \big\uparrow & \big\uparrow & \big\uparrow \\[0.2em]
        \nabla \Theta (x) & \nabla \Sigma_1 (x) & \cdots & \nabla \Sigma_{d-1} (x)\\[0.2em]
        \big\downarrow & \big\downarrow & \big\downarrow & \big\downarrow
    \end{bmatrix} ^\mathsf{T}
    = \Id_{d \times d} \, .
\end{equation}
Moreover, as it is shown in \cite{guillamon2009computational, castejon2013phase, moehliswilsonpre2016}, for points $x = K(\theta, \bs) \in \Omega \subset \mathbb{R}^d$, functions $\nabla \Theta$ and $\nabla \Sigma_i$ satisfy the following adjoint equations
\begin{equation}
    \begin{aligned}\label{eq:prcs_4b}
        \frac{d}{dt} \nabla\Theta(\phi_t(x)) &= -DX^\mathsf{T}(\phi_t(x))\nabla\Theta(\phi_t(x)) \, ,\\
        \frac{d}{dt}\nabla\Sigma_j(\phi_t(x)) &= - \left( DX^\mathsf{T}(\phi_t(x)) - \lambda_j \textrm{Id} \right)\nabla\Sigma_j(\phi_t(x)) \, ,
    \end{aligned}
\end{equation}
which will be used in Section \ref{sec:sec-3} for the computation of the iPRF and iARFs restricted onto the slow submanifold. 

The equations above, when restricted to points on the limit cycle, i.e. $x \in \Gamma$, are the well-known adjoint equations in \cite{ErmentroutKopell91, castejon2013phase}, which have the form 
\begin{equation}\label{eq:NTorder0}
    \frac{1}{T} \frac{d}{d \theta} \nabla \Theta(\gamma (\theta)) = -DX^\mathsf{T}(\gamma(\theta)) \nabla \Theta(\gamma(\theta)) \, ,
\end{equation}
and
\begin{equation}\label{eq:NSorder0}
    \frac{1}{T} \frac{d}{d \theta} \nabla \Sigma_j(\gamma(\theta)) = - \big(DX^\mathsf{T}(\gamma(\theta)) - \lambda_j \Id\big) \nabla \Sigma_j(\gamma(\theta)). 
\end{equation}
Thus, the iPRC and iARC$_j$ are periodic solutions of equations \eqref{eq:NTorder0} and \eqref{eq:NSorder0}, respectively. To guarantee the uniqueness we need to impose the following normalization conditions (see Appendix~\ref{sec:appendix} for more details):
\begin{equation}\label{eq:ncT0}
    \langle \nabla \Theta(\gamma(\theta)), X(\gamma(\theta)) \rangle =\dfrac{1}{T},
\end{equation}
and
\begin{equation}\label{eq:ncS0}
    \langle \nabla \Sigma_j(\gamma(\theta)), K_{e_j}(\theta) \rangle = 1.
\end{equation}

 As with equation \eqref{eq:eqMione}, we can provide the analytical solutions for the adjoint equations \eqref{eq:NTorder0}-\eqref{eq:NSorder0}. Indeed, let $\Psi(t)$ be a fundamental matrix solution of \eqref{eq:NTorder0} and $\Psi_T:=\Psi(T)$ the monodromy matrix of the system. Then, we can obtain periodic solutions of equations \eqref{eq:NTorder0}-\eqref{eq:NSorder0} as long as $e^{-\lambda_j T}$ is an eigenvalue of $\Psi_T$ (see Appendix~\ref{ap:operators}). Therefore, let $(e^{-\lambda_j T}, \hat{w}_j)$ be an eigenpair of the monodromy matrix $\Psi_T$, then a periodic solution of \eqref{eq:NSorder0} is given by 
\begin{equation}\label{eq:fsa}
    \nabla \Sigma_j (\gamma(\theta))=e^{\lambda_j \theta T} \Psi(\theta T) \hat{w}_j.  
\end{equation}
 
The fact that $e^{-\lambda_j T}$ is an eigenvalue of $\Psi_T$ is guaranteed by the relationship between the eigenvalues of $\Phi_T$ defined in \eqref{eq:K1sol_complex} and $\Psi_T$. Indeed, if $e^{\lambda_j T}$ is an eigenvalue of $\Phi_T$ , then $e^{-\overline{\lambda}_j T}$ is an eigenvalue of $\Psi_T$, where $\overline{\lambda}_j$ is the complex conjugate of $\lambda_j$ (see Appendix~\ref{ap:operators} for more details). Therefore, if $\lambda_j \in \mathbb{R}$ then $\overline{\lambda}_j=\lambda_j \in \mathbb{R}$ and $e^{-\lambda_j T}$ is an eigenvalue of $\Psi_T$ whereas, if $\lambda_j \in \mathbb{C}$, then $e^{-\lambda_j T}$ and $e^{-\overline{\lambda}_j T}$ are both eigenvalues of $\Psi_T$. 
In conclusion, if $e^{\lambda_j T}$ is an eigenvalue of the monodromy matrix $\Phi_T$, then $e^{-\lambda_j T}$ is an eigenvalue of the monodromy matrix $\Psi_T$. 

\subsubsection*{Adjoint equations for real amplitude functions}

Adjoint equation \eqref{eq:NSorder0} also holds for real amplitude functions of the form \eqref{eq:icomplexa}, replacing $\lambda_j$ by $\nu_j=\ln |\mu_j|$, that is
\begin{equation}\label{eq:NSorder0_real}
\frac{1}{T} \frac{d}{d \theta} \nabla \widetilde \Sigma_j(\gamma(\theta)) = - \big(DX^\mathsf{T}(\gamma(\theta)) - \nu \Id\big) \nabla \widetilde \Sigma_j(\gamma(\theta)). 
\end{equation}

For the real functions \eqref{eq:real-and-complex-amplitudes}, one can check that the real iARC $\nabla \widetilde \Sigma_j(\gamma(\theta))$ and $\nabla \widetilde \Sigma_{j+1}(\gamma(\theta))$ satisfy the system of equations
\begin{equation}\label{eq:NSorder0_real_case}
    \begin{aligned}
        \frac{1}{T} \frac{d}{d \theta} \nabla \widetilde{\Sigma}_j(\gamma(\theta)) & = -
            DX^\mathsf{T}(\gamma(\theta)) \nabla \widetilde{\Sigma}_j(\gamma(\theta)) + \alpha_{j} \nabla \widetilde{\Sigma}_j(\gamma(\theta)) + \beta_j \nabla \widetilde{\Sigma}_{j+1}(\gamma(\theta)) \, ,
            \\
        \frac{1}{T} \frac{d}{d \theta} \nabla \widetilde{\Sigma}_{j+1}(\gamma(\theta)) & = -
            DX^\mathsf{T}(\gamma(\theta)) \nabla \widetilde{\Sigma}_{j+1}(\gamma(\theta)) - \beta_{j} \nabla \widetilde{\Sigma}_j(\gamma(\theta)) + \alpha_j \nabla \widetilde{\Sigma}_{j+1}(\gamma(\theta)) \, .
    \end{aligned}
\end{equation}

In addition, the real functions $\nabla \widetilde \Sigma_j$ satisfy the normalization conditions
\[\langle \nabla \widetilde \Sigma_j(\gamma(\theta)), \widetilde K_{e_j}(\theta) \rangle =1,\]
where $\widetilde K_{e_j}$ are given in \eqref{eq:Kim} and \eqref{eq:Kpairreal}.

\section{Dynamics restricted to the slow submanifold}\label{sec:sec-3}

In this Section, we present a method to compute the parameterization of the slow submanifold of the periodic orbit, that is, the invariant submanifold associated to the smallest (in modulus) Floquet exponent, $\lambda_{s}:=\lambda_1<0$. For the sake of simplicity, we assume that $\lambda_s$ is real, so the slow manifold is 2-dimensional and orientable. In the complex case, the results also follow but we do not discuss them in detail here (see the Discussion section). Thus, let $\sigma_s=\Sigma_s(x):=\Sigma_1(x) \in \mathbb{R}$ be the amplitude coordinate associated the slow submanifold $\mathcal{S} \subset \mathbb{R}^d$ defined as 
\begin{equation}\label{eq:slow}
    \mathcal{S}:= \{x \in \mathbb{R}^d \enskip | \enskip\Sigma_j(x)=0 \, , \enspace \forall \, j\neq s\}.
\end{equation}

\begin{remark}
Even if we set out our approach for the slowest direction $\sigma_s$, that being associated to the Floquet multiplier closest to 1, our methods will also yield, in a completely analogous way, the parameterization and expressions of the iPRF and iARF restricted to any other direction $\sigma_j$.
\end{remark}

\subsection{Parameterization of the slow submanifold} \label{parameterization_slow_manifold}

The parameterization of the slow manifold is given by 
\begin{equation}\label{eq:K_slow}
    \begin{aligned}
        K_s : \mathbb{T} \times (\widetilde{\mathcal{B}} \subset \mathbb{R}) &\rightarrow \mathbb{R}^d \\
        (\theta, \sigma_s) &\rightarrow K(\theta, \sigma_s,0,\ldots,0) \, .
    \end{aligned}
\end{equation}
Notice that in this case $K(\theta, \sigma_s,0,\ldots,0)$ and $\widetilde K(\theta, \sigma_s,0,\ldots,0)$ coincide. Thus, the map $K_s$ satisfies the invariance equation
\begin{equation}\label{eq:red-InvEq2}
    \frac{1}{T}\frac{\partial}{\partial\theta}K_s(\theta, \sigma_s) + \lambda_s \sigma_s \frac{\partial}{\partial \sigma_s}K_s(\theta, \sigma_s) = X(K_s(\theta, \sigma_s)).
\end{equation}

We address the computation of the parameterization of the slow submanifold $K_s$ analogously to the approach used in the full \( d \)-dimensional case (see Section \ref{sec:computation_K1} and \cite{PerezCervera2020}). We assume that equation \eqref{eq:red-InvEq2} admits a formal series solution of the form 
\begin{equation}\label{eq:red-FourierTaylor}
    K_s(\theta, \sigma_s) = \sum_{n=0}^{\infty} K_n(\theta) \, \sigma_s^{n} \, ,
\end{equation}
and proceed by matching terms with the same powers of $\sigma_s$ on both sides of equation \eqref{eq:red-InvEq2}. This yields a sequence of equations that we solve iteratively, as shown below.

For $n = 0$, the term $K_{\mathbf{0}}(\theta)$ satisfies the equation 
\begin{equation}\label{eq:eqMzero}
    \frac{1}{T}\frac{d}{d\theta}K_{\mathbf{0}}(\theta) = X(K_{\mathbf{0}}(\theta)),
\end{equation}
whose solution is the limit cycle itself, that is $K_{\mathbf{0}}(\theta) = \gamma(\theta T)$.

For $n = 1$, the function $K_1(\theta)$ satisfies the equation
\begin{equation}\label{eq:eqMone}
    \frac{1}{T}\frac{d}{d\theta}K_{1}(\theta) + \lambda_s K_{1}(\theta) = DX(K_{\mathbf{0}}(\theta))K_{1}(\theta),
\end{equation}
whose solution can be expressed as
\begin{equation}\label{eq:eqMzeroSols}
    K_{1}(\theta) = e^{-\lambda_s \theta T } \Phi(\theta T) w_s \, ,
\end{equation}
where $\Phi(t)$ is the solution of the variational equations \eqref{eq:mjVarEqs} and $w_s$ is the eigenvector of the monodromy matrix $\Phi(T)$ associated to the non-trivial slowest decaying Floquet multiplier $\mu_s=e^{\lambda_s T}$.

\begin{remark}
As mentioned in Remark~\ref{rem:nu}, the solutions $K_1(\theta)$ of equation \eqref{eq:eqMone} are not unique. If $w_s$ is an eigenvector of the monodromy matrix, so is $b w_s$ for any $ b \in \mathbb{R}$, which yields new solutions $b K_1 (\theta)$. Even though all the choices of $K_1$ are mathematically equivalent, the choice of $b$ affects the numerical properties of the algorithm. 
\end{remark}

Finally, for $n \geq 2$, the periodic terms $K_{n}(\theta)$ satisfy the so-called \textit{homological equations}
\begin{equation}\label{eq:homologicalEqs}
    \frac{1}{T}\frac{d}{d\theta} K_{n}(\theta) + n \lambda_s K_{n}(\theta) = DX(K_{\mathbf{0}}(\theta)) K_{n}(\theta) + B_{n}(\theta) \, ,
\end{equation}
where $B_{n}(\theta)$ is the coefficient of the term $\sigma_s^n$ in the Taylor expansion of
\begin{equation}\label{eq:Bexpand}
    \emph{X}\left(\sum_{m=0}^{n-1} K_{m}(\theta) \sigma_s^{m} \right).
\end{equation} 
Notice that $B_{n}(\theta)$ is an explicit polynomial depending only on the terms of order lower than $n$, that is, on the functions $K_{m}(\theta)$ for $0\leq m<n$, and whose coefficients are derivatives of $X$ evaluated at $K_{\mathbf{0}}$. They can be numerically computed using automatic differentiation techniques \cite{haro2016}.

The homological equations~\eqref{eq:homologicalEqs} can be solved using formal Fourier series for the functions $K_{n}(\theta)$ where the coefficients of the series are the unknowns. The resulting system of equations for the Fourier coefficients is linear but involves a large-dimensional matrix, leading to a high computational cost when solving it. In the next section we discuss an efficient numerical method to avoid this numerical challenge.

\subsubsection{Reducibility of the homological equations}

In this section, we adopt the strategy proposed in \cite{huguet2013computation} (see also \cite{castelli2015parameterization, PerezSH20})
to solve the homological equations \eqref{eq:homologicalEqs} efficiently in Fourier space. We use a  transformation based on the first order approximation of the parameterization $K$ in \eqref{eq:kThetaSigma} (Floquet normal form) to transform the homological equations \eqref{eq:homologicalEqs} into a diagonal constant-coefficient linear system in Fourier space.

Let us write $K_n$ as
\begin{equation}\label{eq:Kn_form}
    K_n(\theta) = \Q(\theta)\, u(\theta) \, , 
\end{equation} 
where $u(\theta)$ are periodic real vector-valued functions to be determined and $\Q$ the matrix given by
\begin{equation} \label{eq:DK1}
    \Q(\theta) = \begin{bmatrix}
        \big\uparrow & \big\uparrow & \big\uparrow & \big\uparrow \\
        K'_0(\theta) & K_{e_1}(\theta) & \cdots & K_{e_{n-1}}(\theta) \\
        \big\downarrow & \big\downarrow & \big\downarrow & \big\downarrow
    \end{bmatrix} \, ,
\end{equation}
where $K'_0(\theta)$ is the derivative of $K_{\mathbf{0}}(\theta)$ with respect to $\theta$. 
Since the columns of the matrix $\Q(\theta)$ satisfy equations \eqref{eq:eqMizero} and \eqref{eq:eqMione} we have
\begin{equation} \label{eq:edoQQ}
    \frac{1}{T}\frac{d}{d \theta} \Q(\theta) = DX(K_{\mathbf{0}}(\theta)) \Q(\theta) - \Q(\theta) \widetilde \Lambda \, , \quad \text{with} \quad 
    \widetilde \Lambda =\text{diag}(\lambda_0,\lambda_1,\ldots,\lambda_{d-1}).
\end{equation}
where $\lambda_j$ are the Floquet exponents of the periodic orbit (recall $\lambda_0=0$).

Differentiating \eqref{eq:Kn_form} with respect to $t$ and using that $K_n(\theta)$ must satisfy equation \eqref{eq:homologicalEqs} we obtain the following equation
\begin{equation} \label{eq:edouu}
    \frac{1}{T} \left (\frac{d}{d\theta} \Q(\theta) \, u(\theta) + \Q(\theta) \frac{d}{d\theta} u(\theta) \right ) + n \lambda_s \Q(\theta) u(\theta)= DX(K_{\mathbf{0}} (\theta)) \Q(\theta) \, u(\theta) + B_n(\theta) \, .
\end{equation}
Using that matrix $\Q$ satisfies equation \eqref{eq:edoQQ}, equation \eqref{eq:edouu} writes as 
\begin{equation}
      \Q(\theta) \frac{1}{T}\frac{d}{d \theta} u(\theta) - \Q(\theta) \widetilde \Lambda \, u(\theta) = - n \lambda_{s} \Q(\theta) \, u(\theta) + B_n(\theta) \, ,
\end{equation}
and multiplying both sides by $\Q^{-1}(\theta)$ we obtain
\begin{equation} \label{eq:edou_diagonal}
\begin{aligned}
    \frac{1}{T}\frac{d}{d\theta} u(\theta) = ( \widetilde \Lambda - n\lambda_{s} \Id ) \, u(\theta) + A_n(\theta), \\
    \end{aligned}
\end{equation}
where
\begin{equation}
    A_n(\theta) := \mathcal{Q}^{-1}(\theta)B_n(\theta).
\end{equation}

Finally, we write $u(\theta)$ and $A_{n}(\theta)$ in Fourier series, that is,
\begin{equation}\label{eq:mjIneedAgain}
    u(\theta) = \sum_{k = -\infty}^{\infty} u_{k}e^{2\pi ik\theta}, \quad \quad A_{n}(\theta) = \sum_{k = -\infty}^{\infty} A_{k}e^{2\pi ik\theta}, \quad \quad \quad A_{k}, u_{k} \in \mathbb{C}^d,
\end{equation}
and substitute expressions \eqref{eq:mjIneedAgain} into equations \eqref{eq:edou_diagonal}. We obtain a linear system for the Fourier coefficients $u_k=(u_k^{(1)},u_k^{(2)}, \cdots,u_k^{(d)}) \in \mathbb{C}^d$ which is diagonal and can be solved component-wise, thus obtaining the following expression for the Fourier coefficients:
\begin{equation}\label{eq:mjFourierCoefs}
    u^{(j)}_{k} = \frac{1}{\frac{2\pi ik}{T} + n \lambda_s - \lambda_{j-1}}A_{k}^{(j)}, \quad \text{for} \enskip j = 1, 2, \ldots, d \, .
\end{equation}
Notice that the superindex $(j)$ refers to each component of the vectors $u_k$ and $A_k$. Finally, the solution $K_{n}(\theta)$ of equation \eqref{eq:homologicalEqs} is obtained from the transformation \eqref{eq:Kn_form}.

\begin{remark}\label{rem:nonres}
	The Fourier coefficients $u^{(j)}_{k}$ in \eqref{eq:mjFourierCoefs} are formally well defined to all orders provided that, for any $k \in \mathbb{Z}$ and $n\geq 2$, we have
	\[\frac{2\pi ik}{T} + n \lambda_s - \lambda_{j-1} \neq 0, \, \quad j =1, \dots,d.\]
	Notice that this condition is satisfied as long as there are no resonances in the sense of Definition \ref{def:res}, as it implies $\Re(n\lambda_s-\lambda_{j-1}) \neq 0~$ for all $~j=1, \dots,d$.
\end{remark}

\begin{remark} 
    It is possible to design the same strategy using the real functions $\widetilde K_{e_j}$, instead of the complex functions $K_{e_j}$, as columns in the matrix $\mathcal{Q}$. In this case, $\mathcal{Q}$ satisfies equation \eqref{eq:edoQ} with the matrix $\widetilde \Lambda_R:=\text{diag}(\lambda_0,\Lambda_R)$ instead of $\widetilde \Lambda$. Proceeding analogously as for the complex case, it is possible to write the equations for the corresponding real-valued functions $\widetilde u^{(j),(j+1)}$. See Section \eqref{sec:rav} for more details.
\end{remark}

\subsection{Computation of iPRF and iARFs restricted onto the slow submanifold} \label{iPRF_iARF_slow_manifold}
As explained in Section~\ref{sec:phaseAmpFun}, if the full parameterization $x=K(\theta, \bs)$ is known, one can obtain the full set of iPRF and iARFs by means of equation~\eqref{eq:usefulSpace}. However, determining the expressions of the iPRF $\nabla \Theta$ and the iARF $\nabla \Sigma_s$ restricted onto the slow submanifold $\mathcal{S}$ in \eqref{eq:slow} parameterized by $K_s$ in \eqref{eq:K_slow} requires a different approach. In this section, we will compute them by means of solving the adjoint equations formally. Thus, we assume a power series solution in $\sigma_s$
\begin{equation}\label{eq:nT_expand} 
    \begin{aligned}
        \nabla \Theta (K_s(\theta,\sigma_s)) &= \sum^\infty_{n=0} Z_n(\theta)\sigma_s^n, \\
        \nabla \Sigma_s (K_s(\theta,\sigma_s)) &= \sum^\infty_{n=0} I_n(\theta)\sigma_s^n \, ,
    \end{aligned}
\end{equation}
for the equations
\begin{equation}\label{eq:adjoint}
    \begin{aligned}
        \frac{d}{dt} \nabla \Theta (K_s(\theta,\sigma_s))&=-DX^\mathsf{T}(K_s(\theta,\sigma_s)) \nabla \Theta (K_s(\theta,\sigma_s)), \\
        \frac{d}{dt} \nabla \Sigma_s (K_s(\theta,\sigma_s))&= -[DX^\mathsf{T}(K_s(\theta,\sigma_s)) - \lambda_s \Id] \nabla \Sigma_s (K_s(\theta,\sigma_s)) \, ,
    \end{aligned}
\end{equation}
obtained from the adjoint equations \eqref{eq:prcs_4b} restricted to points $x=K_s(\theta,\sigma_s)$ on the slow manifold.

Using automatic differentiation \cite{haro2016}, we can compute expansions in $\sigma_s$ of $DX^\mathsf{T}(K(\theta,\sigma_s))$, that is, 
\begin{equation}\label{eq:ad_exp}
    DX^\mathsf{T}(K_s(\theta,\sigma_s)):= \sum^\infty_{n=0} F_n(\theta) \sigma_s^n ,
\end{equation}
where $F_n(\theta)$ are $d \times d$ matrices whose expressions depend on the terms $K_m(\theta)$ for $m=0, \ldots, n$. Notice that $F_0(\theta)=DX^\mathsf{T}(K_{\mathbf{0}}(\theta))$.

Substituting the power series \eqref{eq:nT_expand} into equations \eqref{eq:adjoint}, equating terms of the same order in $\sigma_s$, and using that $\dot{\sigma}_s = \lambda_s \sigma_s$, we obtain the homological equations for $\nabla \Theta$ and $\nabla \Sigma_s$ up to any order $n \geq 0$: 

\begin{equation}\label{eq:NTordern}
\frac{1}{T} \frac{d}{d \theta} Z_n(\theta) = -\big[DX^\mathsf{T}(K_{\mathbf{0}}(\theta)) + n \lambda_s \Id \big] Z_n(\theta) - G_n(\theta) \, ,
\end{equation}

\begin{equation}\label{eq:NSordern}
\frac{1}{T} \frac{d}{d \theta} I_n(\theta) = -[ DX^\mathsf{T}(K_{\mathbf{0}}(\theta)) + (n-1)\lambda_s \Id] I_n(\theta) - H_n(\theta),
\end{equation}
where
\begin{equation*}
    G_n(\theta)=\sum_{i = 0}^{n-1} F_{n-i}(\theta) Z_i(\theta)
    \quad \text{and} \quad H_n(\theta)= \sum_{i=0}^{n-1} F_{n-i}(\theta) I_{i} (\theta) \, , \text{ for } n \geq 1
\end{equation*}
and $G_0(\theta) \equiv 0$ and $H_0(\theta) \equiv 0$.

\begin{remark}
    Notice that for $n=0$ equations \eqref{eq:NTordern} and \eqref{eq:NSordern} correspond to the adjoint equations \eqref{eq:NTorder0} and \eqref{eq:NSorder0}, respectively, since $Z_0(\theta)=\nabla \Theta (\gamma(\theta))$ and $I_0(\theta)=\nabla \Sigma_s (\gamma(\theta))$. The equations can be solved uniquely when the normalization conditions \eqref{eq:ncT0} and \eqref{eq:ncS0} are imposed. 
\end{remark}

Notice that equation \eqref{eq:NSordern} for $n=1$ is of the form 
\begin{equation}\label{eq:NTorder1}
\frac{1}{T} \frac{d}{d \theta} I_1(\theta) = -DX^\mathsf{T}(K_{\mathbf{0}}(\theta)) I_1(\theta) - F_1(\theta)I_0(\theta) \, .
\end{equation}
Thus, it is a periodic linear non-homogeneous equation, whose homogeneous part corresponds to the adjoint equation \eqref{eq:NTorder0} which admits periodic solutions of the form 
\[I^h_1(\theta)=c Z_0(\theta) \, , \textrm{ for } c \in \mathbb{R} \, .\]
Therefore, periodic solutions of equation \eqref{eq:NTorder1} can be written as 
\[I_1(\theta)=c Z_0(\theta) + I^p_1(\theta) \, ,\]
where $I^p_1(\theta)$ is a particular periodic solution of \eqref{eq:NTorder1} and $c$ is an arbitrary constant. Hence, periodic solutions are not unique. For this reason we need to add a normalization condition (see Appendix~\ref{sec:appendix}) given by
\begin{equation}\label{eq:Nc-Wilson-v2}
    \langle I_0(\theta),K_1'(\theta)\rangle + \langle I_1(\theta),X(K_{\mathbf{0}}(\theta))\rangle=0 \,,
\end{equation}
or equivalently, using equations \eqref{eq:ncT0} and \eqref{eq:eqMone}, 
\begin{equation}\label{eq:Nc-Wilson}
    \langle I_0(\theta),DX(K_{\mathbf{0}}(\theta))K_1(\theta) \rangle + \langle I_1(\theta),X(K_{\mathbf{0}}(\theta)) \rangle =\lambda_s \, ,  
\end{equation}
which is the expression considered in \cite{Wilson2020}.

Next, we proceed to solve equations \eqref{eq:NTordern} and \eqref{eq:NSordern} efficiently for $n\geq1$ following a similar strategy to the one presented in \cite{PerezSH20} to solve the equations for $K_n$.

\subsubsection{Efficient solution of the adjoint equations}

In this section we derive an efficient procedure for computing the $n^{\text{th}}$-order terms, $Z_n$ and $I_n$, of the iPRF and the iARF expanded along the least-contracting direction, associated to $\lambda_{\text{s}}$ (see equation \eqref{eq:nT_expand}). 

The functions $Z_n(\theta)$ and $I_n(\theta)$ can be expressed as a linear combination of the $0^{\text{th}}$-order terms, $\nabla \Theta^0 (\theta) :=\nabla \Theta(\gamma(\theta))$ and $\nabla \Sigma_j^0 (\theta):=\nabla \Sigma_j(\gamma(\theta))$, for $j=1, \ldots, d-1$ (see Appendix \ref{ap:operators}). Thus, we write
\begin{equation}\label{eq:les-adjuntes}
    Z_n(\theta) = Q(\theta)\, v(\theta) \, , \quad \text{and} \quad I_n(\theta) = Q(\theta)\, w(\theta) \, , 
\end{equation} 
where $v(\theta)$ and $w(\theta)$ are scalar periodic real-valued functions to be determined and $Q(\theta)$ the $d \times d$ matrix with the solutions of the adjoint equations \eqref{eq:NTorder0} and \eqref{eq:NSorder0}, arranged as
\begin{equation}\label{eq:matrixQ}
    Q(\theta) := \begin{bmatrix}
        \big\uparrow & \big\uparrow & \big\uparrow & \big\uparrow \\[0.4em]
        \nabla \Theta^0(\theta) & \nabla \Sigma_1^0(\theta) & \cdots & \nabla \Sigma_{d-1}^0(\theta) \\[0.4em]
        \big\downarrow & \big\downarrow & \big\downarrow & \big\downarrow
    \end{bmatrix} \, .
\end{equation}
Therefore, the matrix $Q(\theta)$ satisfies the differential equation
\begin{equation} \label{eq:edoQ}
    \frac{1}{T}\frac{d}{d \theta} Q(\theta) = - DX^\mathsf{T}(K_{\mathbf{0}}(\theta)) Q(\theta) + Q(\theta) \widetilde \Lambda \, ,
\end{equation}
where $\widetilde \Lambda$ is the diagonal matrix in \eqref{eq:edoQQ}. The case with real functions $ \nabla \widetilde\Sigma^0_j$ is discussed in Section \ref{sec:rav}.

Differentiating \eqref{eq:les-adjuntes} with respect to $t$ (recall that $\theta=t/T$) and using that $Z_n(\theta)$ and $I_n(\theta)$ satisfy equations \eqref{eq:NTordern} and \eqref{eq:NSordern}, respectively, we obtain 
\begin{equation} \label{eq:edow}
    \begin{aligned}
        \frac{1}{T} \left (\frac{d}{d\theta} Q(\theta) \, v(\theta) + Q(\theta) \frac{d}{d\theta} v(\theta) \right )&= - \bigl[DX^\mathsf{T}(K_{\mathbf{0}} (\theta)) + n\lambda_{\text{s}} \Id \bigr] Q(\theta) \, v(\theta) - G_n(\theta), \\
        \frac{1}{T} \left ( \frac{d}{d\theta} Q(\theta) \, w(\theta) + Q(\theta) \frac{d}{d\theta} w(\theta) \right )&= - [ DX^\mathsf{T}(K_{\mathbf{0}}(\theta)) +(n-1) \lambda_{s} \Id ] Q(\theta) \, w(\theta) - H_n(\theta) \, .
    \end{aligned}
\end{equation}
Using that $Q(\theta)$ fulfills equation \eqref{eq:edoQ}, system \eqref{eq:edow} reduces to
\begin{equation}
    \begin{aligned}
         Q(\theta) \widetilde \Lambda \, v(\theta) + Q(\theta) \frac{1}{T}\frac{d}{d \theta} v(\theta) &= -  n \lambda_{s} Q(\theta) \, v(\theta) - G_n(\theta) \, ,\\
          Q(\theta) \widetilde \Lambda \, w(\theta) + Q(\theta) \frac{1}{T}\frac{d}{d\theta} w(\theta) &= -  (n-1) \lambda_{s} Q(\theta) \, w(\theta) -H_n(\theta) \, ,
    \end{aligned}
\end{equation}
and multiplying both sides of the equations by $Q^{-1}(\theta)$ we obtain that $v(\theta)$ and $w(\theta)$ satisfy the following equations:
\begin{equation} \label{eq:edow_diagonal}
    \begin{aligned}
        \frac{1}{T}\frac{d}{d\theta} v(\theta) &= - \bigl[\widetilde \Lambda + n\lambda_{s} \Id\bigr] \, v(\theta) + \mathcal{G}_n(\theta), \\
        \frac{1}{T} \frac{d}{d\theta} w(\theta) &= 
        - \bigl[\widetilde \Lambda + (n-1)\lambda_{s} \Id\bigr] \, w(\theta) + \mathcal{H}_n(\theta) \, ,
    \end{aligned}
\end{equation}
where
\begin{equation}
    \mathcal{G}_n(\theta) := -Q^{-1}(\theta)G_n(\theta), \qquad \mathcal{H}_n(\theta) := -Q^{-1}(\theta)H_n(\theta) \, .
\end{equation}
Recall that matrix $Q(\theta)$ is invertible because its column vectors constitute a basis of the space of solutions of the adjoint equations (see Appendix \ref{ap:operators}).

\begin{remark}
Notice that it is not necessary to invert matrix $Q$ if the normal bundle to the periodic orbit is known. Indeed, using equation \eqref{eq:usefulSpace} (see also Appendix C in \cite{PerezCervera2020}) we have that 
\[Q^{-1}(\theta)=\mathcal{Q}^\mathsf{T}(\theta) \, ,\]
where $\mathcal{Q}$ is given in \eqref{eq:DK1}.
\end{remark}

Finally, we use Fourier series to solve equations \eqref{eq:edow_diagonal}. Indeed, we express functions $v, w, \mathcal{G}_n$ and $\mathcal{H}_n$ in Fourier series as
\begin{equation}\label{eq:fourier-dec}
    \begin{aligned}
        v(\theta) = \sum_{k = -\infty}^\infty v_k e^{2 \pi i k \theta} \, ,& \qquad 
        \mathcal{G}_n(\theta) = \sum_{k = -\infty}^\infty g_k e^{2 \pi i k \theta} \, , \\
        w(\theta) = \sum_{k = -\infty}^\infty w_k e^{2 \pi i k \theta} \, ,& \qquad
        \mathcal{H}_n(\theta) = \sum_{k = -\infty}^\infty h_k e^{2 \pi i k \theta} \, , \qquad v_k, w_k, g_k, h_k \in \mathbb{C}^d
    \end{aligned}
\end{equation}
and substitute these expressions back in \eqref{eq:edow_diagonal}. We obtain a linear system for the Fourier coefficients
$v_k, w_k \in \mathbb{C}^d$, which is diagonal and can be solved component-wise, thus obtaining
the following expression for the Fourier coefficients:
\begin{equation}\label{eq:sol-z1}
    v^{(j)}_k = \frac{g^{(j)}_k}{ \frac{2\pi i k}{T} + \lambda_{j-1} + n\lambda_{s} }, 
\end{equation}
and
\begin{equation}\label{eq:sol-i1}
    w^{(j)}_k = \frac{h^{(j)}_k}{  \frac{2\pi i k}{T} + \lambda_{j-1} + (n-1) \lambda_{s} }, \enspace \text{~for} \enskip j=1,\ldots , d.
\end{equation}
Here, the superscript $(j)$ refers to each component of the vectors $v_k$, $w_k$, $g_k$ and $h_k$. Finally, once the functions $v$ and $w$ have been determined for order $n$, we compute $Z_n$ and $I_n$ for $n \geq 1$ in a straightforward way as indicated by expressions \eqref{eq:les-adjuntes}.

\begin{remark}\label{rem:fcoef1}
    The Fourier coefficients $v_k^{(j)}$ in \eqref{eq:sol-z1} are formally well defined for $n \geq 1$ provided that for any $k \in \mathbb{Z}$ and $n \geq 1$, we have
    \[2\pi i \frac{k}{T} + \lambda_{j-1} + n\lambda_{s} \neq 0, \enskip \text{for} \enskip j=1,\ldots d. \]
    This condition is always satisfied since we assumed that $\Re(\lambda_j) < 0$ with the exception of $\lambda_0=0$, but in the latter case, the denominator does not vanish either because $\lambda_s \neq \lambda_0$ and $n \geq 1$. 
\end{remark}

\begin{remark}\label{rem:fcoef2}
The Fourier coefficients $w_k^{(j)}$ in \eqref{eq:sol-i1} are formally well defined for $n \geq 1$ provided that for any $k \in \mathbb{Z}$ and $n \geq 1$, we have
\[\frac{2\pi i  k}{T} + \lambda_{j-1} + (n-1)\lambda_{s} \neq 0, \enskip \text{for} \enskip j=1,\ldots d. \]
The condition is not satisfied when $k=0$, $j=1$ ($\lambda_{0}=0$) and $n=1$ (see next remark). For the rest of the cases the condition is always satisfied, since we assumed that $\Re(\lambda_j) < 0$ (with the exception of $\lambda_0=0$, but $\lambda_s \neq \lambda_0$). 
\end{remark}

\begin{remark}
When $k=0$, $j=1$ ($\lambda_{0}=0$) and $n=1$ the denominator in \eqref{eq:sol-i1} vanishes. In this case, a solution exists if $h_0^{(1)}=0$ and, therefore, $w_0^{(1)}$ is free and can be chosen according to the normalization condition \eqref{eq:Nc-Wilson}. 
\end{remark}

\subsubsection{Real amplitude variables}\label{sec:rav}

It is possible to design the same strategy using the real functions $\widetilde \Sigma_{j}$ defined in \eqref{eq:icomplexa} and \eqref{eq:real-and-complex-amplitudes} instead of the complex ones $\Sigma_{j}$. Thus, we consider a real matrix $\widetilde Q$, whose columns are the real-valued functions $\nabla \widetilde{\Sigma}_{j}^{0} (\theta):=\nabla \widetilde{\Sigma}_{j}(\gamma(\theta))$, for the change of coordinates \eqref{eq:les-adjuntes}. Then, using that $\nabla \widetilde{\Sigma}^0_{j}$ satisfies equation \eqref{eq:NSorder0_real} when $\lambda_j=\nu + i \pi/T$ and $\nabla \widetilde{\Sigma}^0_{j,j+1}$ satisfy system \eqref{eq:NSorder0_real_case} when $\lambda_{j,j+1}$ are complex conjugates, we have that the real matrix $\widetilde Q$ satisfies system 
\begin{equation} \label{eq:edoQtilde}
    \frac{1}{T}\frac{d}{d \theta} \widetilde Q(\theta) = - DX^\mathsf{T}(K_{\mathbf{0}}(\theta)) \widetilde Q(\theta) + \widetilde Q(\theta) \widetilde \Lambda_R \, ,
\end{equation}
where $\widetilde \Lambda_R := \text{diag} (\lambda_0, \Lambda_R^{\mathsf{T}})$, being $\Lambda_R$ the matrix defined in \eqref{eq:Lambda_real}.

Proceeding analogously as for the complex case, the Fourier coefficients for the corresponding real-valued functions $\widetilde v^{(j),(j+1)}$ are determined by solving the $2\times2$ linear system
\begin{equation}\label{eq:Fcoef_real}
    \begin{pmatrix}
       \widetilde v_k^{(j)} \\[0.3em]
       \widetilde v_k^{(j+1)}
    \end{pmatrix}
    =
    \begin{pmatrix}
        2 \pi i \frac{k}{T} + \alpha_j + n \lambda_s & -\beta_j \\[0.3em]
        \beta_j & 2 \pi i \frac{k}{T} + \alpha_j + n \lambda_s
    \end{pmatrix}^{-1}
    \begin{pmatrix}
       \widetilde g_k^{(j)} \\[0.3em]
       \widetilde g_k^{(j+1)}
    \end{pmatrix}
    =:
    M^{-1}
    \begin{pmatrix}
       \widetilde g_k^{(j)} \\[0.3em]
       \widetilde g_k^{(j+1)}
    \end{pmatrix}
    \, .
\end{equation}
\begin{remark}
Notice that system \eqref{eq:Fcoef_real} can be solved uniquely as long as $\det M \neq 0$. Indeed, we have
\[\det{M}=\alpha^2+\beta^2 + 2 \alpha (2 \pi i \frac{k}{T} +n \lambda_s) + (2 \pi i \frac{k}{T} + n \lambda_s)^2. \]
Letting $\xi:=2 \pi i \frac{k}{T} + n \lambda_s$ and using that $\lambda_j \overline \lambda_{j}=\alpha_j^2+\beta_j^2$ and $\lambda_j + \overline \lambda_j=2\alpha_j$, we can write
\[\det{M}=\lambda_j \overline \lambda_j + \xi (\lambda_j + \overline \lambda_j) + \xi^2 = (\xi + \lambda_j)(\xi+ \overline \lambda_j) \neq 0.\]
Since $\lambda_{j+1}=\overline{\lambda}_j$, this is equivalent to the conditions in Remark~\ref{rem:fcoef1}.
\end{remark}

Analogously, a similar expression can be found for the corresponding Fourier coefficients of $\widetilde w^{(j),(j+1)}$. Namely,
\begin{equation}\label{eq:Fcoef_real_w}
    \begin{pmatrix}
       \widetilde w_k^{(j)} \\[0.3em]
       \widetilde w_k^{(j+1)}
    \end{pmatrix}
    =
    \begin{pmatrix}
        2 \pi i \frac{k}{T} + \alpha_j + (n-1) \lambda_s & -\beta_j \\[0.3em]
        \beta_j & 2 \pi i \frac{k}{T} + \alpha_j + (n-1) \lambda_s
    \end{pmatrix}^{-1}
    \begin{pmatrix}
       \widetilde h_k^{(j)} \\[0.3em]
       \widetilde h_k^{(j+1)}
    \end{pmatrix} \, .
\end{equation}

\begin{remark}
    Similarly to the previous remark, system \eqref{eq:Fcoef_real_w} can be solved as long as the conditions in Remark \ref{rem:fcoef2} are satisfied.
\end{remark}

\section{Numerical application}\label{sec:numerics}

In this section, we present the results of the numerical implementation of the methodology described in Section~\ref{sec:sec-3}, applied to a relevant model in computational neuroscience.
The model considered is a 6-dimensional system of differential equations, introduced in \cite{dumont2019macroscopic}, that describes the mean-field dynamics of a neural network consisting of excitatory and inhibitory cells (E-I network). For a given set of the parameters, the system presents a hyperbolic limit cycle, whose Floquet multipliers take real positive, negative and complex values. Thus, the system represents a comprehensive example of most of the remarkable key points of the methodology. 

The system of differential equations for the model in \cite{dumont2019macroscopic} is given by
\begin{equation} \label{eq:dumont_gutkin_model}
    \begin{aligned}
        \tau_e \frac{d}{d t} r_e &= \frac{\Delta_e}{\pi \tau_e} + 2 r_e V_e \, , \\
        \tau_e \frac{d}{d t} V_e &= V_e^2 + \eta_e - (\tau_e \pi^2 r_e)^2 - \tau_e S_{ei} + I_e^{\text{ext}} \, , \\
        \tau_{si} \frac{d}{d t} S_{ei} &= -S_{ei} + J_{ei} r_i \, , \\
        \tau_i \frac{d}{d t} r_i &= \frac{\Delta_i}{\pi \tau_i} + 2 r_i V_i \, , \\
        \tau_i \frac{d}{d t} V_i &= V_i^2 + \eta_i - (\tau_i \pi^2 r_i)^2 + \tau_i S_{ie} + I_i^{\text{ext}} \, , \\
        \tau_{se} \frac{d}{d t} S_{ie} &= -S_{ie} + J_{ie} r_e \, . \\
    \end{aligned}
\end{equation}
For the set of parameters
\begin{equation*}
    \mathcal{P} = \{\tau_{e,i} = 10,\, \Delta_{e,i} = 1,\, \eta_{e,i} = -5,\, \tau_{se,si} = 1,\, J_{ei}=J_{ie} = 15,\, I_e^{\text{ext}} = 10,\, I_i^{\text{ext}} = 0\} \, ,
\end{equation*}
the system has a hyperbolic attracting limit cycle $\Gamma$. The Floquet multipliers $\mu_j$ and the corresponding Floquet exponents $\lambda_j$ (unique except for the sum of multiples of $2 \pi  \mathrm{i}/T$) associated to $\Gamma$ are
\begin{equation*}
    \begin{aligned}
        \boldsymbol{\mu_0} &= 1 \,, \\
        \boldsymbol{\mu_1} &= 0.0537 \,, \\
        \boldsymbol{\mu_2} &= 2.3\cdot10^{-4} + 3.13\cdot10^{-4} \,\mathrm{i} \,, \\
        \boldsymbol{\mu_3} &= 2.3\cdot10^{-4} - 3.13\cdot10^{-4} \,\mathrm{i} \,, \\
        \boldsymbol{\mu_4} &= -3.99\cdot10^{-10} \,, \\
        \boldsymbol{\mu_5} &= -1.57\cdot10^{-10} \,,
    \end{aligned}
    \hspace{5em}
    \begin{aligned}
        \boldsymbol{\lambda_0} &= 0 \,, \\
        \boldsymbol{\lambda_1} &= -0.1405 \,, \\
        \boldsymbol{\lambda_2} &= -0.377 + 0.045 \,\mathrm{i} \,, \\
        \boldsymbol{\lambda_3} &= -0.377 - 0.045 \,\mathrm{i} \,, \\
        \boldsymbol{\lambda_4} &= -1.040 + \mathrm{i} \pi/T \,  \,, \\
        \boldsymbol{\lambda_5} &= -1.0845 + \mathrm{i} \pi/T \, \,.
    \end{aligned}
\end{equation*}

Notice that the Floquet multipliers $\mu_1 \in \mathbb{R}^+$, $\mu_2$ and $\mu_3$ are a pair of complex conjugated eigenvalues and $\mu_4, \mu_5 \in \mathbb{R}^{-}$, thus covering all the cases discussed in Section~\ref{sec:computation_K1}. Notice also that there are no resonances in the Floquet exponents.

\begin{figure}[htbp!]
    \centering
    \begin{tabular}{cc}
        \includegraphics[width=0.5\textwidth]{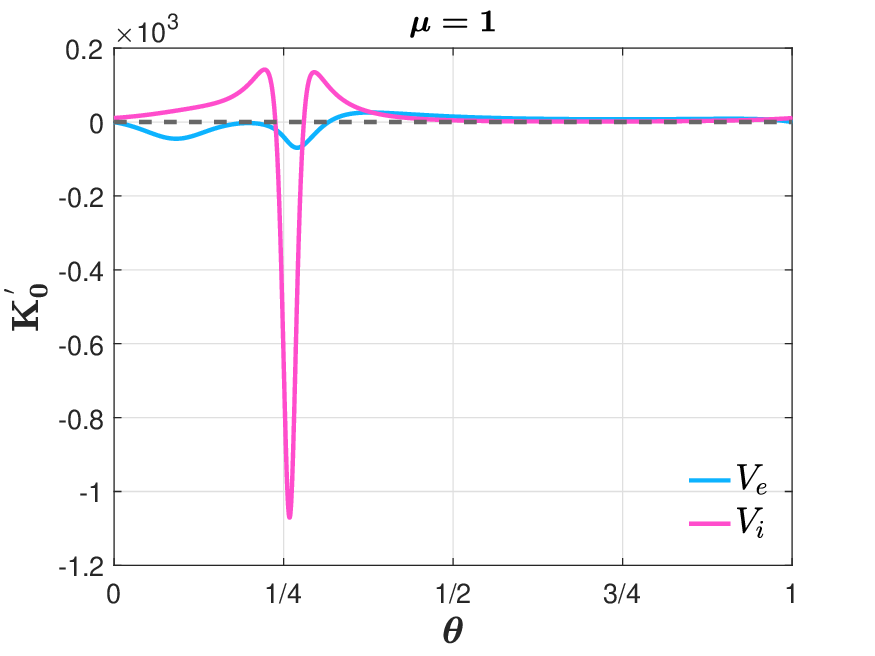} & \includegraphics[width=0.5\textwidth]{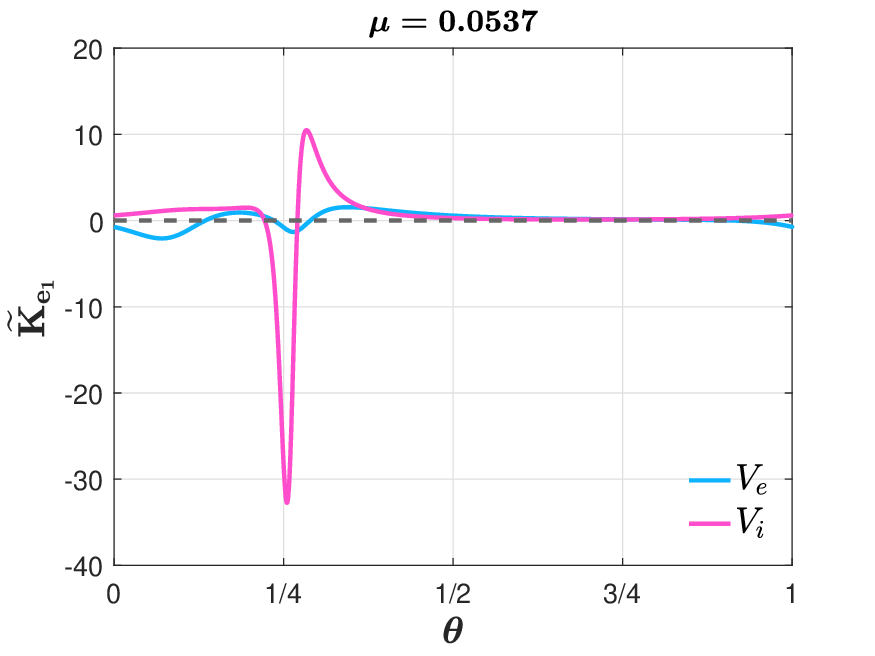} \\
        \includegraphics[width=0.5\textwidth]{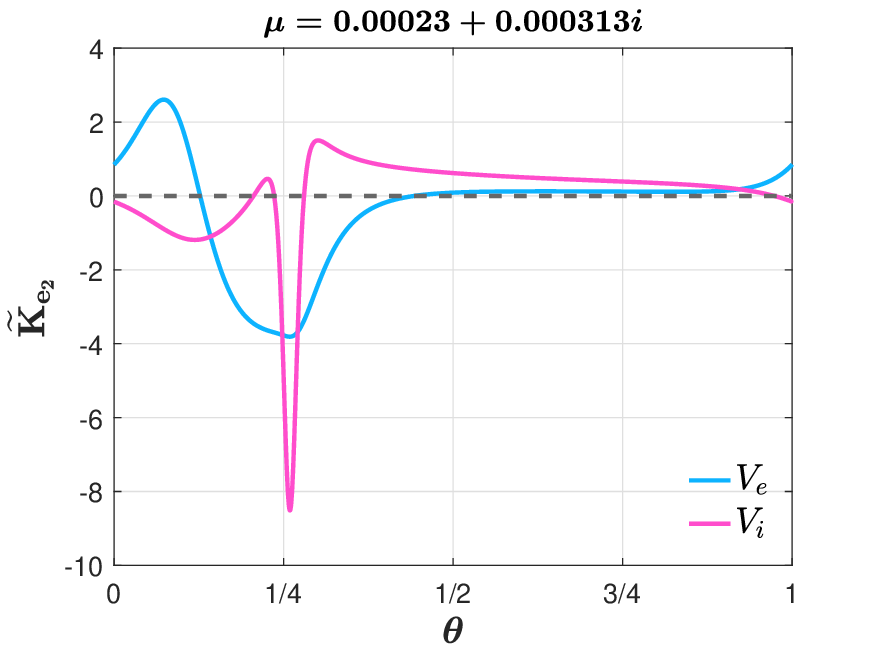} & \includegraphics[width=0.5\textwidth]{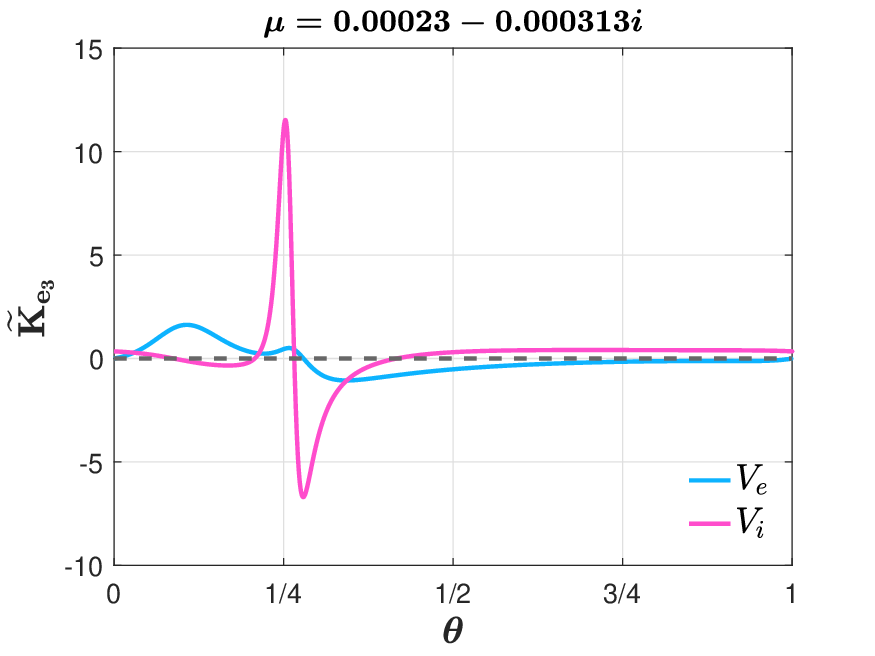} \\
        \includegraphics[width=0.5\textwidth]{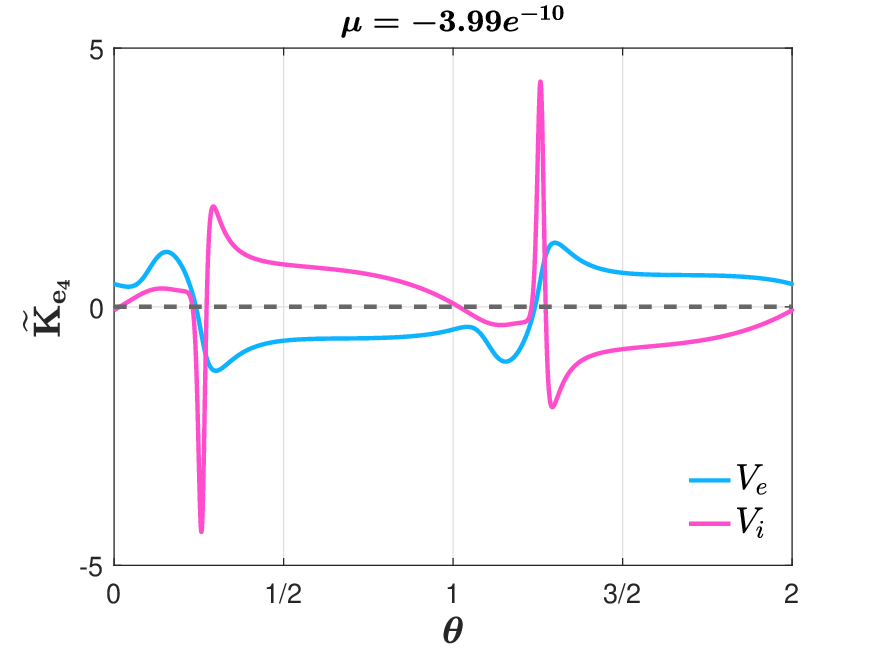} & \includegraphics[width=0.5\textwidth]{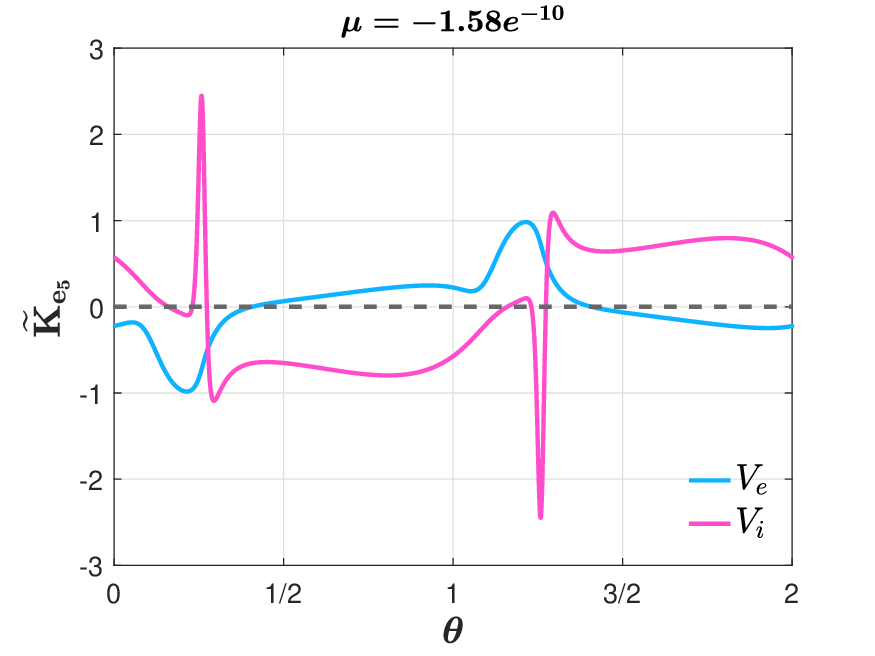}
    \end{tabular}
    \caption{Tangent and normal bundles of the periodic orbit $\Gamma$ of system \eqref{eq:dumont_gutkin_model}. We plot the coordinates $V_e$ (blue) and $V_i$ (cyan) of the real functions $K_{\mathbf{0}}^{\prime}$ (tangent bundle) and $\widetilde K_{e_j}$, $j=1, \ldots, 5$, (normal bundle), obtained from formulas \eqref{eq:K1sol_complex} for real positive multipliers, \eqref{eq:Kim} for real negative multipliers and \eqref{eq:Kpairreal} for complex conjugate ones. Notice that real functions $\widetilde K_{e_4, e_5}$ associated to negative Floquet multipliers $\mu_{4,5}$ are 2-periodic.}
    \label{fig:K1_solutions}
\end{figure}

We implemented the methodology to compute a local approximation of the parameterization of the slow submanifold associated with the least-contracting Floquet multiplier $\mu_1$ as well as the iPRFs and iARFs restricted to it. Thus, we first computed the tangent and normal bundle associated to the periodic orbit, which constitute the change of coordinates in \eqref{eq:DK1} to solve the equations for $K_n$ efficiently (see Section~\ref{parameterization_slow_manifold}). In Figure~\ref{fig:K1_solutions}, we show the $V_e$ and $V_i$ components of the periodic real-valued functions $K_{\mathbf{0}}^{\prime}(\theta)$ and $\widetilde K_{e_j}(\theta)$, $1 \leq j \leq 5$, of the formal series $\eqref{eq:Kformal_real}$. 

\begin{figure}[htbp!]
    \centering
    \begin{tabular}{cc}
        \includegraphics[width=0.5\textwidth]{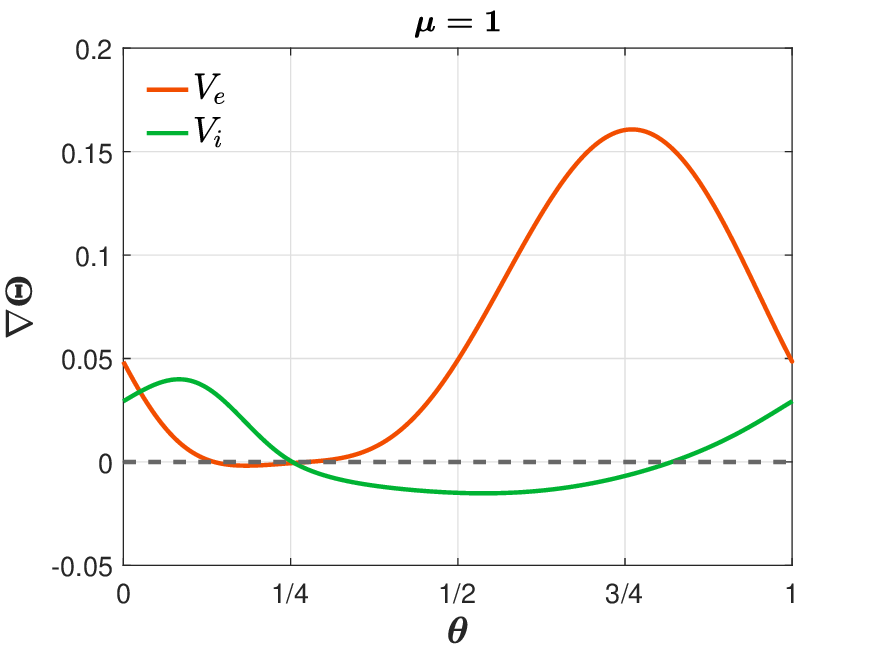} & \includegraphics[width=0.5\textwidth]{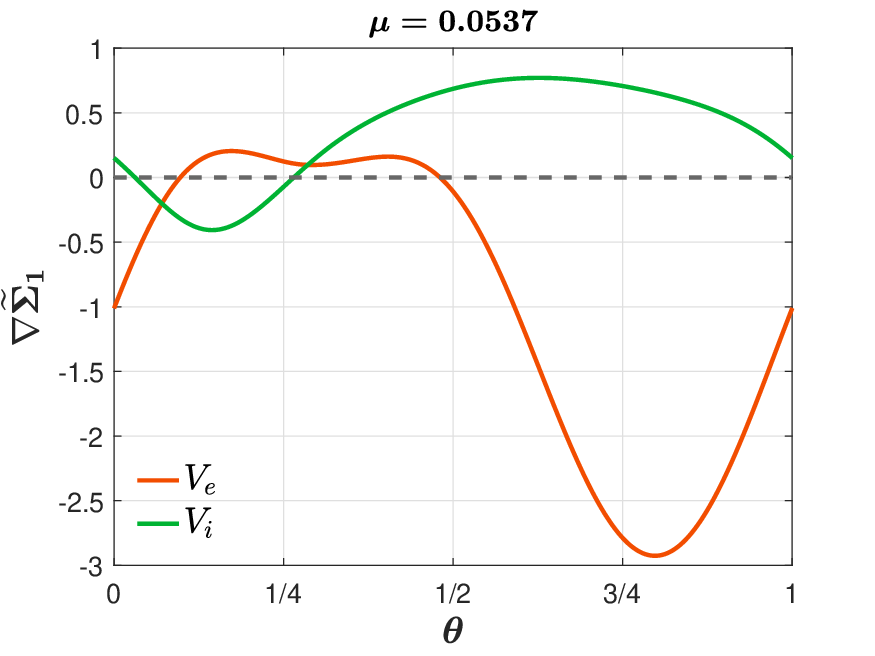} \\
        \includegraphics[width=0.5\textwidth]{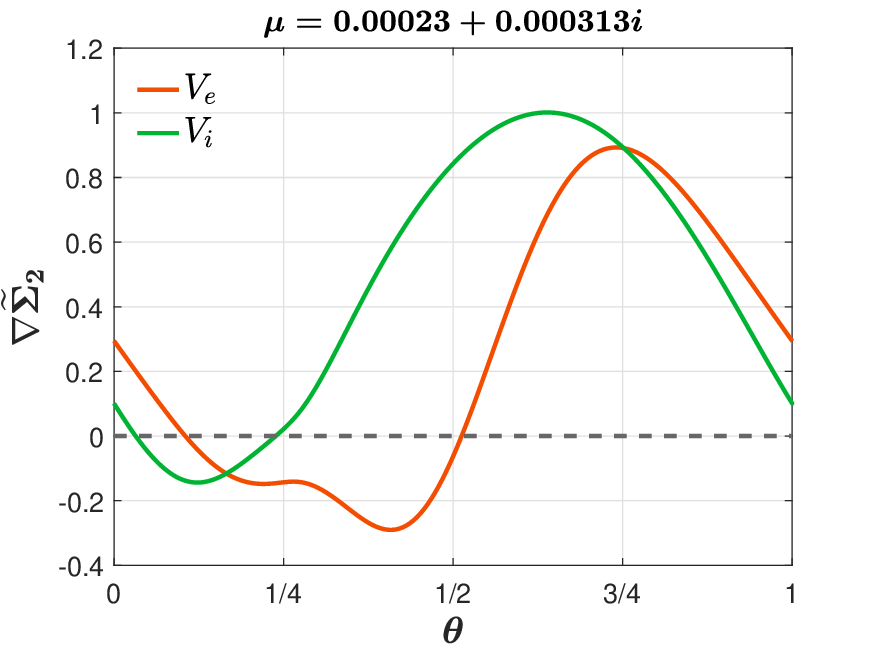} & \includegraphics[width=0.5\textwidth]{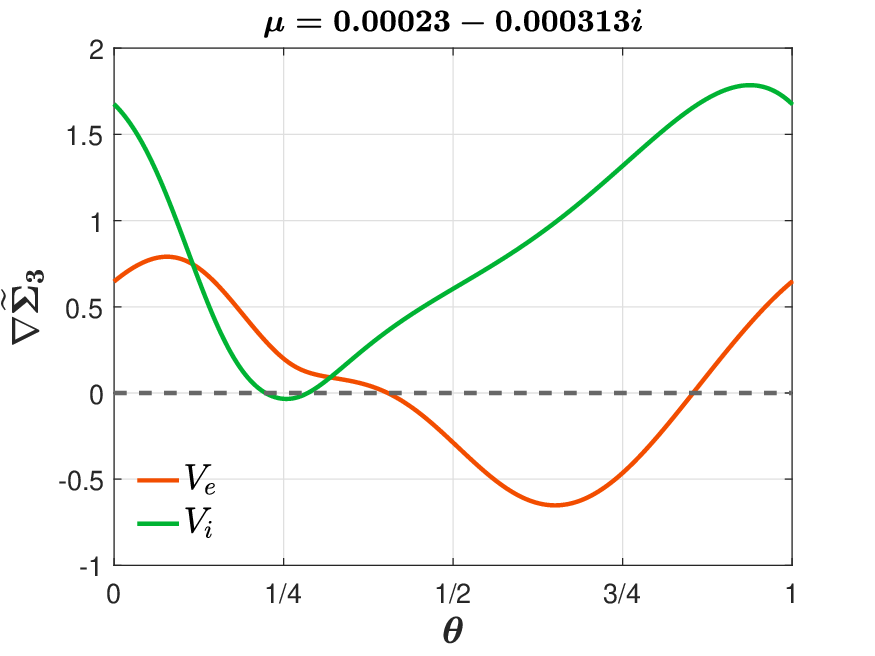} \\
        \includegraphics[width=0.5\textwidth]{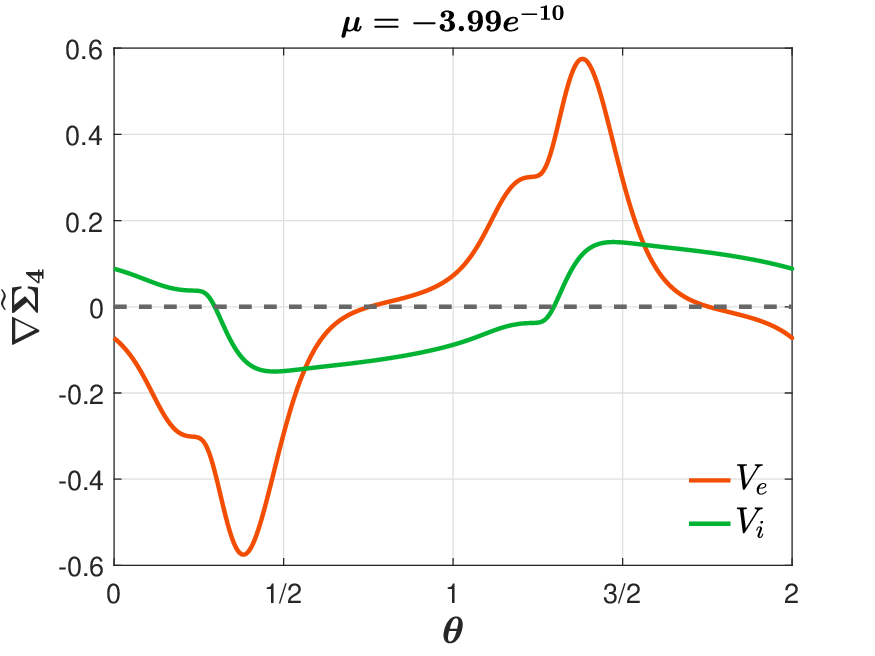} & \includegraphics[width=0.5\textwidth]{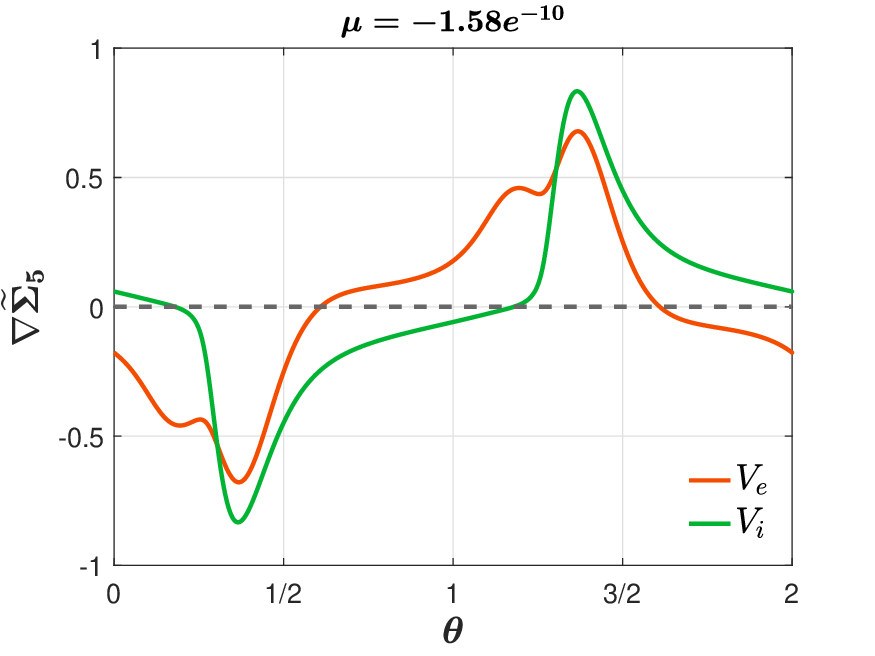}
    \end{tabular}
    \caption{Functions iPRF and iARFs along the periodic orbit $\Gamma$ of system \eqref{eq:dumont_gutkin_model}. We plot the coordinates $V_e$ (orange) and $V_i$ (green) for the iPRC $\nabla \Theta(\gamma(\theta))$ and the iPRCs $\nabla \widetilde \Sigma_j(\gamma(\theta))$, $j=1, \ldots, 5$, the latter being obtained from expressions \eqref{eq:fsa} for real positive Floquet multipliers, \eqref{eq:NSorder0_real} for real negative multiplers and \eqref{eq:NSorder0_real_case} for complex conjugate ones. Notice that real functions $\nabla \widetilde \Sigma_{4,5}$ associated to negative Floquet multipliers $\mu_{4,5}$ are 2-periodic.}
    \label{fig:nabla0}
\end{figure}

Another important set of functions for our methodology, see Section \ref{sec:phaseAmpFun}, is the $0^{\text{th}}$-order terms of the iPRF and iARFs, that is, the classical iPRC $\nabla \Theta(\gamma(\theta))$ and iARCs $ \nabla \widetilde \Sigma_j(\gamma(\theta))$, for $1 \leq j \leq 5$, as they define the change of coordinates \eqref{eq:les-adjuntes} needed to solve the equations for $Z_n$ and $I_n$ efficiently. We compute them by Taylor expanding expression \eqref{eq:usefulSpace} in $\bs$ and solving the $0^{\text{th}}$-order equation (as in \cite{PerezCervera2020}). In Figure \ref{fig:nabla0} we show the $V_e$ and $V_i$ components of the computed iPRC and iARCs.

\begin{figure}[htbp!]
    \centering
    \begin{tabular}{ll}
        \large{\textbf{A}} & \large{\textbf{B}} \\
        \includegraphics[width=0.48\textwidth]{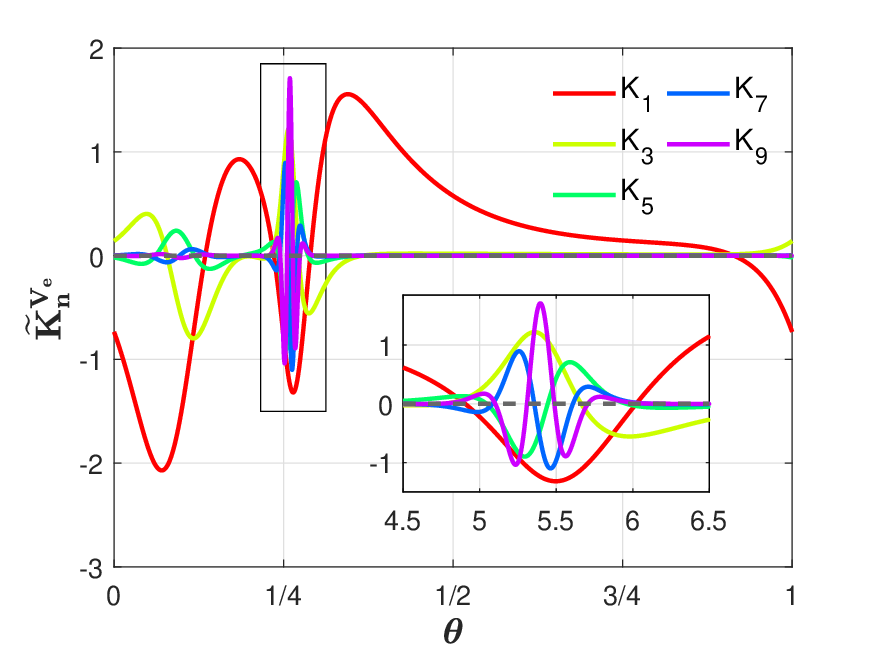} & \includegraphics[width=0.5\textwidth]{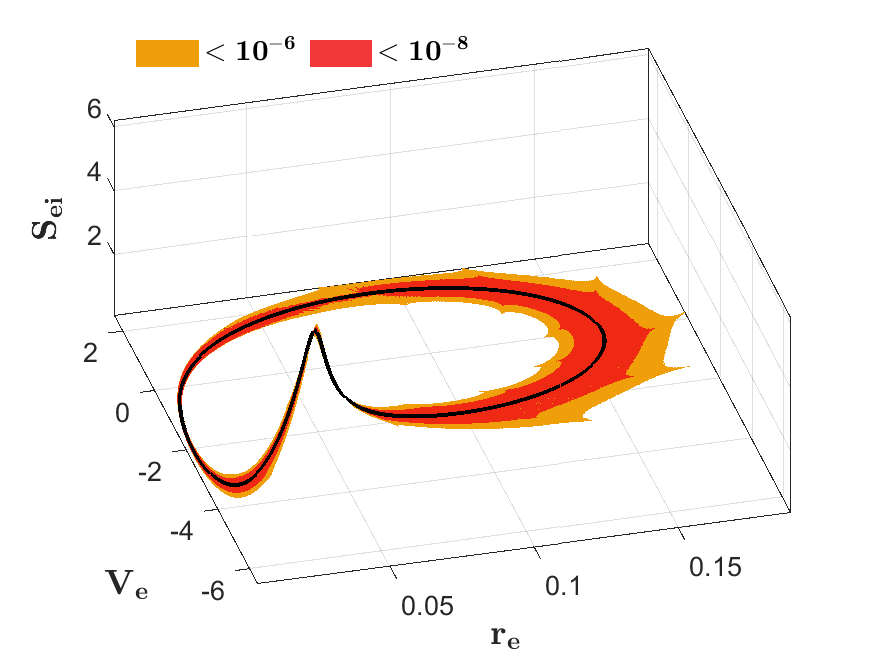} \\
    \end{tabular}
    \caption{Graphical overview of the numerical results for the slow manifold's local parameterization for system~\eqref{eq:dumont_gutkin_model}. (A) Coordinate $V_e$ of the coefficient functions $K_n(\theta)$, for $n$ odd, of the local parameterization $K_{L,N}(\theta,\sigma)$ of \eqref{eq:red-FourierTaylor}, truncated at order $L=9$ with $N=2^{12}$ (see colour legend). Notice how these functions decay to zero as the order increases, except around phase $1/4$ where functions vary abruptly (see inset plot). (B) Projection of the slow manifold's local parameterization onto the $(r_e,V_e,S_{e i})$-coordinated system for two domains of accuracy: one with a tolerance $10^{-8}$ (red surface) and another with a smaller tolerance of $10^{-6}$ (orange surface). The limit cycle is depicted by a solid-black curve.}
    \label{fig:Kn_coordVe}
\end{figure}

We compute a local approximation of the parameterization of the slow submanifold \eqref{eq:red-FourierTaylor} using high-order Fourier-Taylor expansions. To demonstrate the efficiency of the strategy presented in Section~\ref{parameterization_slow_manifold}, we computed power expansions up to order $L=9$ and Fourier expansions with $N=2^{12}$ coefficients to approximate periodic functions $K_n$ (see Figure~\ref{fig:Kn_coordVe}A). Thus, the truncated function 
\[K_{L,N}(\theta,\sigma):=\sum_{n=0}^{L} \sum_{k=-N/2}^{N/2} C_k \, e^{2 \pi \textrm{i} k \theta} \, \sigma^n, \]
provides an approximation of the parameterization $K_s(\theta,\sigma)$ of the slow submanifold in a neighbourhood of the periodic orbit. To determine its domain of accuracy we define the error function associated to $K_{L,N}$ as
\begin{equation}\label{error_inv_eq}
    E_{L,N}(\theta, \sigma) = \left\lVert \, \sum_{n=0}^{L} \bigg[ \frac{1}{T} K_n^{\prime}(\theta) + n \lambda_s K_n(\theta) \bigg]\sigma^n - X\bigg( \sum_{n=0}^{L} K_n(\theta) \sigma^n \bigg) \, \right\rVert_2 \, .
\end{equation}
where $\lVert \cdot \rVert_{L^2}$ denotes the $L_2$-norm.

The derivatives of functions $K_n$ are calculated via Fast Fourier Transform (FFT) with as many Fourier coefficients as in the computation of $K_n$. The domain of accuracy is determined as
\[\Omega_{tol}=\{(\theta,\sigma) \in \mathbb{T} \times \mathbb{R} \enskip | \enskip E_{L,N}(\theta, \sigma) < E_{tol}\},\] 
for a preset tolerance $E_{tol}$. In Figure \ref{fig:Kn_coordVe}, we display the $V_e$-component of the functions $K_n$ for some orders (see Figure~\ref{fig:Kn_coordVe}A) as well as the projection of the local approximation $K_{L,N}$ of the slow manifold onto the $(r_e,V_e,S_{ei})$-space, for two domains of accuracy corresponding to  two different tolerances $E_{tol}=10^{-6}$ (orange surface) and $E_{tol}=10^{-8}$ (red surface) (see Figure \ref{fig:Kn_coordVe}B).

\begin{figure}[htbp!]
    \centering
    \begin{tabular}{ll}
        \large{\textbf{A}} & \large{\textbf{B}} \\
        \includegraphics[width=0.48\textwidth]{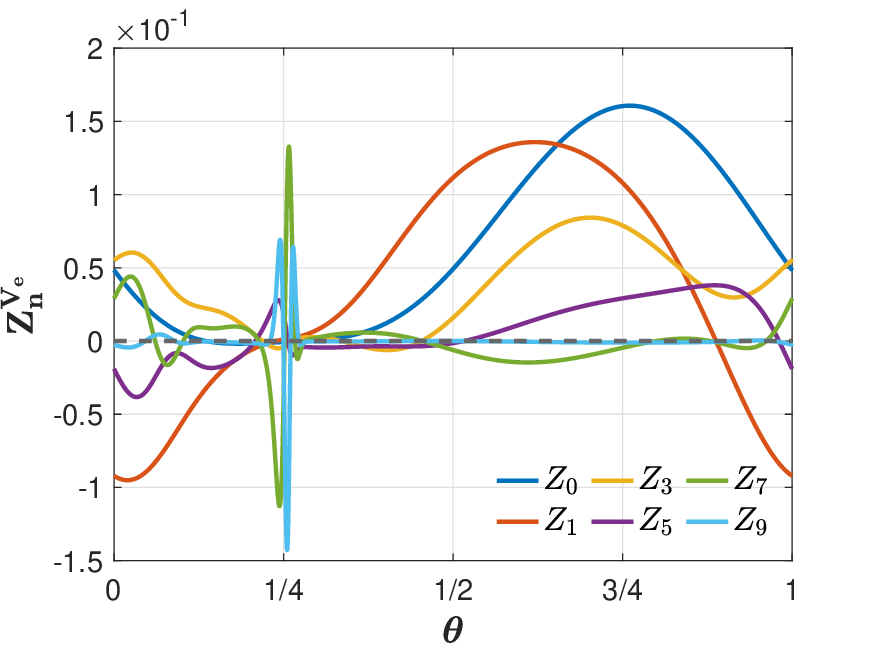} & \includegraphics[width=0.5\textwidth]{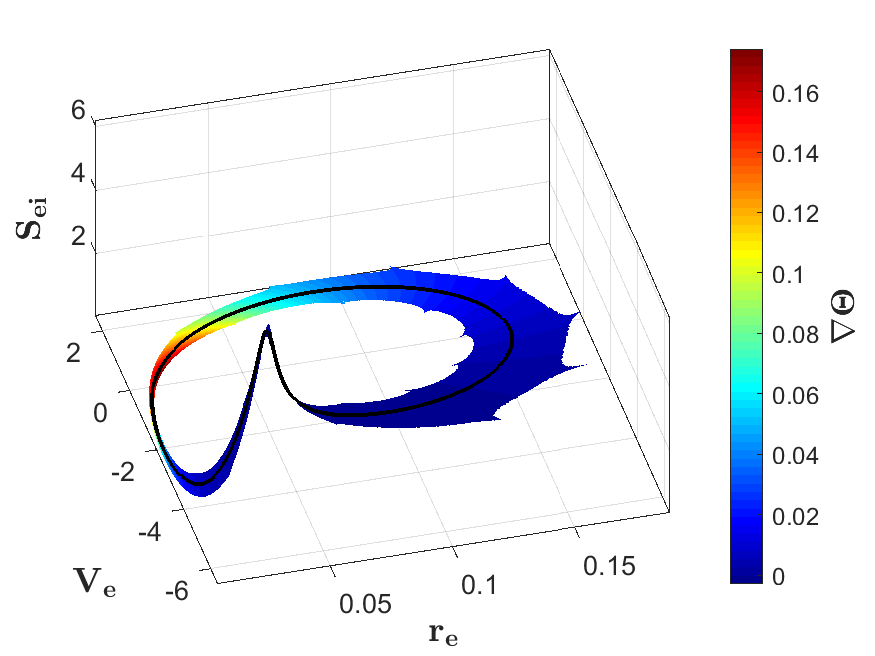} \\ 
        \large{\textbf{C}} & \large{\textbf{D}} \\
        \includegraphics[width=0.48\textwidth]{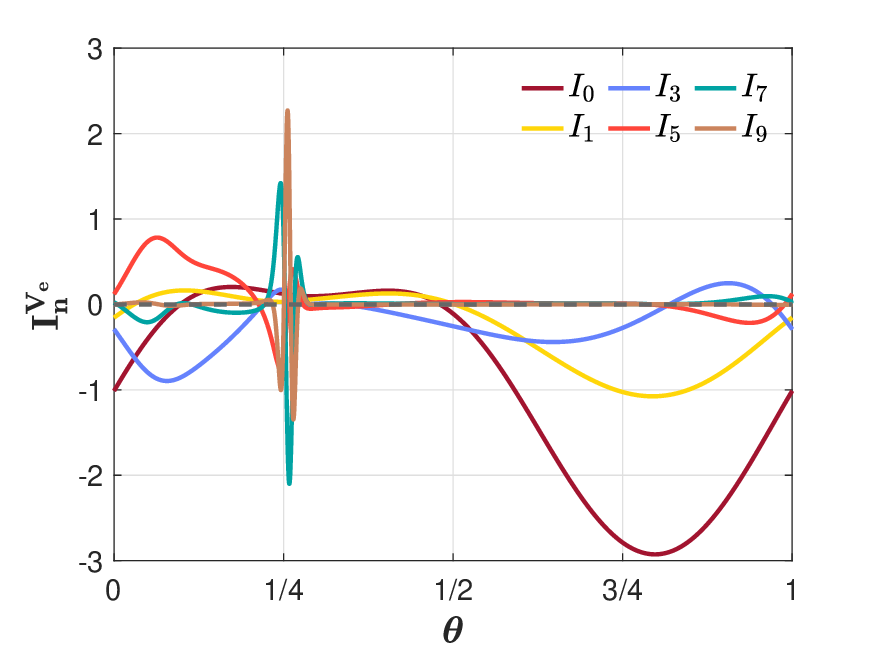} & \includegraphics[width=0.5\textwidth]{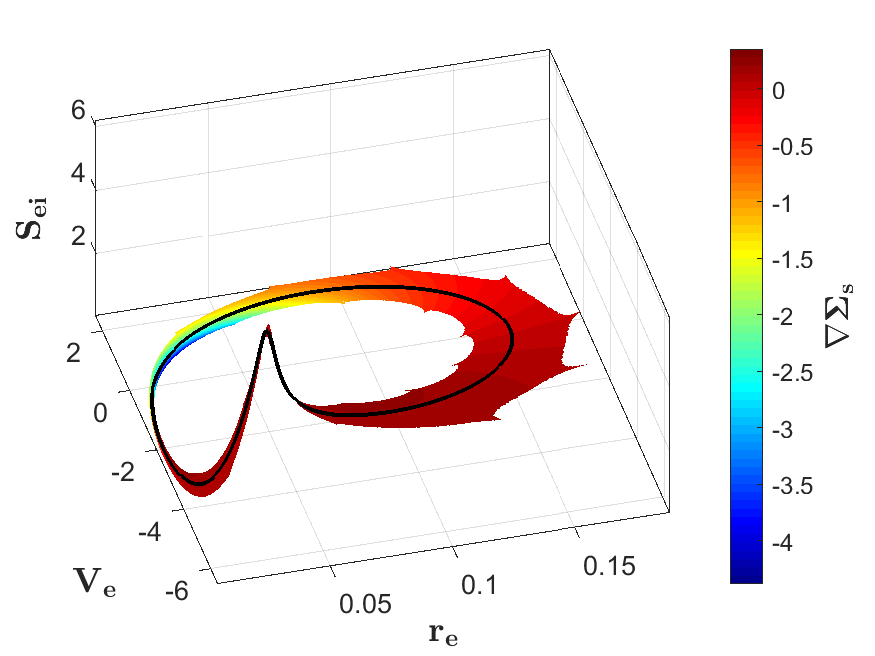}
    \end{tabular}
    \caption{Graphical overview of the numerical results for the approximated iPRF and iARF on the local approximation of the slow submanifold for system~\eqref{eq:dumont_gutkin_model}. (Left column) Coordinate $V_e$ of the functions $Z_n$ (A) and $I_n$ (C) for $n$ odd (see legends). Functions $Z_n$, for $n = 1,3,5,7,9$, have been scaled each by factors $\kappa=10, 100, 10^3, 10^4, 10^4$ while functions $I_n$, for $n=3,5,7,9$, have been scaled by $\kappa=25,250, 10^3, 10^3$, respectively. (Right column) Projections of the slow manifold's local approximation onto the $(r_e,V_e,S_{ei})$-coordinated system, coloured according to the $V_e$-component of the vector-valued functions iPRF (B) and the iARF (D) on it. The black curve represents the limit cycle.}
    \label{fig:Zn_In_solutions}
\end{figure}

Likewise, we also computed an equally accurate approximation (same order trunction in the Fourier-Taylor series) of the iPRF and iARFs (see equations \eqref{eq:nT_expand}) on the locally approximated slow manifold (which corresponds to orange surface in Figure \ref{fig:Kn_coordVe}B). In Figure \ref{fig:Zn_In_solutions} we show the $V_e$-component of the $n^{\text{th}}$-order terms $Z_n$ and $I_n$ in the Taylor expansions of the iPRF $\nabla \Theta$ and the iARF $\nabla \Sigma_s$ (panels A and C, respectively) alongside the $V_e$ component of each one for the points on the local approximation of the slow submanifold (see colour palette in right column). In Figure~\ref{fig:Zn_In_solutions}A and C, the terms $Z_n$ and $I_n$ of higher orders have been scaled by different factors so that they are all of comparable magnitude. As before, the functions $Z_n$ and $I_n$ have been approximated by Fourier series with $N = 2^{12}$ coefficients.

Regions of the manifold where the iPRF (resp. iARF) takes larger values correspond to regions where the level curves of $\Theta$ (isochrons) (resp. $\Sigma$ (isostables)) are more densely packed in space. In such regions, the gradient of the respective function is larger, meaning that any infinitesimal perturbation applied along the $V_e$ direction at a point $(\theta,\sigma)$ has a greater effect on phase or amplitude compared to regions where the level curves are more widely spaced. For the oscillator described by system~\eqref{eq:dumont_gutkin_model}, both phase and amplitude exhibit strong sensitivity to perturbations in the $V_e$-direction in regions where the approximated manifold narrows around the limit cycle (see Figure~\ref{fig:Zn_In_solutions} B and D). 
Furthermore, the foliation of isostables reveals that contraction to the limit cycle along the slow manifold is not uniform along the isochrons.

Notice that the differential equations \eqref{eq:NTordern} and \eqref{eq:NSordern} for $Z_n$ and $I_n$, respectively, involve the non-homogeneous terms $G_n$ and $H_n$. Computing these functions requires, in turn, to compute the $n^{\text{th}}$-order matrices $F_n$ of the Taylor expansion in $\sigma$  of matrix $DX^\mathsf{T}(K_s(\theta,\sigma))$ (see expansion \eqref{eq:ad_exp}). Typically, matrices $F_n$ are obtained by means of automatic differentiation techniques (see \cite{haro2016}). However, for this particular model, this expansion can be derived analytically at any order, as all the entries of matrix $DX^\mathsf{T}$ are linear in its variables. Therefore, the analytical expression for the $n^{\text{th}}$-order terms $F_n$, $n \geq 0$, is given straightaway by 
\[
F_n(\theta)=
\begin{pmatrix}
    \dfrac{2}{\tau_e} K^{V_e}_n(\theta)& -2 \tau_e \pi^2 K^{r_e}_n(\theta) & 0 & 0 & 0 & \delta[n]\dfrac{J_{ie}}{\tau_{se}} \\[0.9em]
    \dfrac{2}{\tau_e} K^{r_e}_n(\theta) & \dfrac{2}{\tau_e} K^{V_e}_n(\theta) & 0 & 0 & 0 & 0 \\[0.9em]
    0 & - \delta[n] & \delta[n] (- \dfrac{1}{\tau_{si}}) & 0 & 0 & 0 \\[0.9em]
    0 & 0 & \delta[n]\dfrac{J_{ei}}{\tau_{si}} & \dfrac{2}{\tau_i} K^{V_i}_n(\theta) & -2 \tau_i \pi^2 K^{r_i}_n(\theta) & 0 \\[0.9em]
    0 & 0 & 0 & \dfrac{2}{\tau_i} K^{r_i}_n(\theta) & \dfrac{2}{\tau_i} K^{V_i}_n(\theta) & 0 \\[0.9em]
    0 & 0 & 0 & 0 & \delta[n] & \delta[n] (-\dfrac{1}{\tau_{se}}) \\
\end{pmatrix} \, ,
\]
where $\delta[n]=1$ if $n=0$ and $0$ otherwise.

\section{Discussion}\label{sec:discussion}

In this paper, we have introduced an efficient numerical method for computing a local parameterization of the slow submanifold (when it exists) of a hyperbolic attracting limit cycle in high-dimensional systems. In addition, we also described the perturbed dynamics restricted to this submanifold in terms of phase and amplitude coordinates, specifically through the infinitesimal Phase Response Function (iPRF) and Amplitude Response Functions (iARFs).
Our methodology builds upon prior results in solving functional equations via coordinate transformations, which enable the development of efficient numerical algorithms \cite{huguet2013computation, PerezCervera2020, James16, Huguet2012, haro2016}.

In earlier work \cite{PerezCervera2020}, we computed the parameterization $K(\theta,\bs)$ of the whole (stable) invariant manifold of limit cycles.  The iPRF and the full set of iARFs were obtained directly by inverting the Jacobian matrix $DK$. By setting $\sigma_j$, for $j=2,\ldots,d-1$, to zero in $K(\theta,\bs)$, the slow submanifold and its dynamics were clearly identified. However, this approach is computationally expensive for high-dimensional systems  when the focus is on a specific direction. Indeed, computing the full parameterization in such cases requires determining all $K_\alpha$ of \eqref{eq:formal} for every combination of $\alpha$, resulting in unnecessary computational overhead. 

To overcome this issue, we formulate an invariance equation for the slow submanifold. Our approach provides an effective and efficient method that only requires knowledge of the normal bundle in the remaining directions. Furthermore, we extend this strategy beyond the invariant manifold itself, providing efficient methods to obtain solutions of the adjoint equations that yield the restricted iPRF $\nabla \Theta$ and iARF $\nabla \Sigma$ on the slow submanifold.

This method offers significant computational advantages over alternative approaches \cite{Wilson2020, Wilson2019}, which rely on direct integration of the differential equations \eqref{eq:NTordern} and \eqref{eq:NSordern}. By choosing an appropriate coordinate system, we demonstrate that the homological equations \eqref{eq:homologicalEqs} reduce to diagonal, constant-coefficient systems in Fourier space, avoiding the inversion of dense matrices, a process that is both computationally expensive and numerically unstable (see \cite{guillamon2009computational}). Building on this idea, we introduce a change of coordinates (in this case provided by the iPRC and iARC) to similarly simplify the adjoint equations. This reduces their solution to a linear system that remains diagonal in Fourier space, ensuring efficiency and stability.

Beyond its performance benefits, this method is particularly advantageous for its broad applicability, especially in the context of high-dimensional systems. The implementation requires only a few key components: the expansion of the differential matrix $DX$ restricted to the (sub)manifold under study, the normal bundle of the periodic orbit, and the iPRC and iARCs.

Our study focuses on the slow manifold to capture the transient dynamics of trajectories as they approach the limit cycle. However, the methodology is not limited to this specific manifold. Indeed, it is also applicable to other invariant (sub)manifolds corresponding to the remaining amplitude coordinates $\sigma_j$, $j = 2,\ldots,d-1$. 

We emphasize that our results are based on the assumption that the Floquet exponents of the limit cycle are non-resonant, as defined in Definition \ref{def:res}. While the presence of a resonance in the Floquet exponents at some order does not prevent the existence of a stable invariant manifold or, consequently, the isochrons (see, for instance, \cite{guckenheimer2013nonlinear}), it does prevent the existence of an analytic conjugacy that transforms the vector field $X$ into a linear one \eqref{eq:phase-amplitude}. In such cases, the parameterization method remains applicable, but the conjugacy must be established with a polynomial vector field rather than a linear one, as done in this study. For a detailed discussion of resonances and the parameterization method in the context of fixed points and equilibria, we refer the reader to \cite{cabre2003parameterization}.

Our methodology does not assume the Floquet exponents are real; it also includes the cases where the Floquet exponents are complex. In this case, one can work with either complex or real amplitude variables. In this study, we have addressed both cases, presenting the relationship between real and complex amplitude formulations (see also \cite{Wilson2019}).

For simplicity, we have considered the slow submanifold to be associated with a real eigenvalue, making it 2-dimensional and orientable (see Section \ref{sec:sec-3}). However, our methodology readily extends to more general cases. For instance, it applies to situations where the largest (in modulus) Floquet multiplier is negative, resulting in a non-orientable bundle. It also applies to cases involving a pair of complex conjugate eigenvalues, where the manifold is 3-dimensional, and the dynamics on the linear bundle is rotational. In the latter case, the parameterization of the manifold requires two complex amplitude variables.

In this paper, we have provided the computation of local approximations of the slow submanifold and the functions evaluated on it but we have not discussed the globalization of them via backwards numerical integration (see \cite{PerezCervera2020, huguet2013computation}). This is a delicate problem, particularly when there are strongly contracting directions, since backwards integration amplifies small errors in the direction of the dominant eigenvalue (see \cite{simo1990analytical}). 
This makes globalization through backwards integration to be numerically unstable. In \cite{simo1990analytical, PerezCervera2020} one can find some strategies to overcome these problems or use a different approach based on splines \cite{osinga2010continuation, Krauskopfetal05}. 

In recent years, there has been growing interest in reconstructing the dynamics of coupled oscillator systems from data. Initially, such reconstructions relied on phase reduction methods \cite{Rosenblum2001, Yeldesbay2019, Kralemann2014}, but more recently, data-driven phase-amplitude reconstruction techniques have been developed \cite{williams2015data, Wilson2020b, kaiser2021data, Wilson2023, Cestnik2022, Yeldesbay2024}. We hope that the methods presented in this paper contribute to advancing this field by providing more effective tools for such reconstructions.

In conclusion, this paper introduces efficient numerical methods to restrict oscillatory dynamics to the slow attracting invariant submanifold. By parameterizing this manifold and computing the infinitesimal Phase and Amplitude Response Functions (iPRF and iARF), we provide a practical framework for simplifying and accurately describing complex oscillatory systems. We expect that the efficient techniques presented here will offer valuable insights and alternative perspectives for studying high dimensional oscillatory dynamics across a wide range of applications, including data-driven phase-amplitude approaches.

\subsection*{Acknowledgements}
Work produced with the support of the grant PID-2021-122954NB-I00 funded by MCIN/AEI/ 10.13039/501100011033 and “ERDF: A way of making Europe”, the Maria de Maeztu Award for Centers and Units of Excellence in R\&D (CEX2020-001084-M) and the AGAUR project 2021SGR1039. The author APC is a Serra Húnter Fellow. 
We also acknowledge the use of the cluster of the UPC Dynamical Systems group for research computing \url{https://dynamicalsystems.upc.edu/en/computing/}.

\appendix

\section{Orthogonality relationships between terms of the asymptotic expansions}\label{sec:appendix}

In this section we present the orthogonality relationships satisfied by the terms of the asymptotic expansions of the parameterization of the slow manifold $K_s$ in \eqref{eq:red-FourierTaylor} and the functions $\nabla \Theta$ and $\nabla \Sigma_s$ evaluated on the slow manifold (see equations \eqref{eq:nT_expand}). These relationships will be used as the normalization conditions \eqref{eq:ncT0}, \eqref{eq:ncS0} and \eqref{eq:Nc-Wilson} to solve solve equations \eqref{eq:NTorder0}, \eqref{eq:NSorder0} and \eqref{eq:NTorder1}, respectively, uniquely. This relationships can be derived analogously for the parameterization $K$ satisfying equation \eqref{eq:mjInvEq} (see also \cite{Wilson2020}).

From now on, we omit the subscript $s$ in $K$ and $\Sigma$. Taking derivatives on both sides of expressions 
\begin{equation}\label{eq:Teq}
\Theta (K(\theta,\sigma))=\theta,
\end{equation}
and
\begin{equation}\label{eq:Seq}
\Sigma (K(\theta,\sigma))=\sigma,
\end{equation}
we have
\begin{equation}\label{eq:NTd}
\langle \nabla \Theta (K(\theta,\sigma)), \frac{\partial K}{\partial \theta} (\theta,\sigma) \rangle =1, \qquad 
\langle \nabla \Theta (K(\theta,\sigma)), \frac{\partial K}{\partial \sigma} (\theta,\sigma) \rangle =0,  
\end{equation}
and
\begin{equation}\label{eq:NSd}
\langle \nabla \Sigma (K(\theta,\sigma)), \frac{\partial K}{\partial \theta} (\theta,\sigma) \rangle =0, \qquad 
\langle \nabla \Sigma (K(\theta,\sigma)), \frac{\partial K}{\partial \sigma} (\theta,\sigma) \rangle =1,
\end{equation}
where $\langle \cdot,\cdot \rangle $ indicates the dot product. Using that
\[ \frac{\partial}{\partial \theta} K (\theta,\sigma)= \sum_{n=0}^{\infty} K'_n(\theta) \sigma^n, \qquad 
\frac{\partial}{\partial \sigma} K (\theta,\sigma)= \sum_{n=0}^{\infty} K_n(\theta) n \sigma^{n-1},\]
(here $'$ denotes derivative with respect to $\theta$) and
equating order by order, from expressions \eqref{eq:NTd} we obtain
\begin{equation}\label{eq:ncZ1}
\langle Z_0,K'_0(\theta) \rangle =1, \qquad \sum_{i = 0}^{n} \langle Z_i(\theta), K'_{n-i}(\theta) \rangle =0, \textrm{ for } n \geq 1
\end{equation}
and
\begin{equation}\label{eq:ncZ2}
\sum_{i = 0}^{n} \langle Z_i(\theta), (n+1-i)K_{n+1-i}(\theta) \rangle =0, \textrm{for } n \geq 0,
\end{equation}
and from expressions \eqref{eq:NSd}, we obtain,
\begin{equation}\label{eq:ncI1}
\sum_{i = 0}^{n} \langle I_i(\theta), K'_{n-i}(\theta) \rangle =0, \textrm{ for } n \geq 0,
\end{equation}
and
\begin{equation}\label{eq:ncI2}
\langle I_0(\theta),K_1(\theta) \rangle =1, \qquad \sum_{i=0}^n \langle I_i(\theta), (n+1-i)K_{n+1-i}(\theta) \rangle =0, \textrm{ for } n \geq 1.
\end{equation}

Notice that expressions~\eqref{eq:ncZ1} and \eqref{eq:ncI2} (for $n=0$) correspond to the normalization conditions \eqref{eq:ncT0} and \eqref{eq:ncS0}, respectively, used to solve equations \eqref{eq:NTorder0} and \eqref{eq:NSorder0} uniquely. Moreover, expression \eqref{eq:ncI1} for $n=1$ corresponds to the normalization condition \eqref{eq:Nc-Wilson-v2} used to solve equation~\eqref{eq:NTorder1} uniquely. The other orthogonality relationships are already satisfied by solving the corresponding equations, which have a unique solution. See the discussion in Section~\ref{ap:operators}.

\section{Solutions of the Adjoint Equations \eqref{eq:NTorder0}-\eqref{eq:NSorder0}}\label{ap:operators}

In this section we discuss how to obtain a periodic solution satisfying the adjoint equations \eqref{eq:NTorder0}-\eqref{eq:NSorder0}. Following \cite{cabre2005parameterization} we introduce the operator 
\[\Lo^*:= \frac{1}{T} \frac{d}{d \theta} + DX^\mathsf{T}(K_{\mathbf{0}} (\theta)). \]
From Proposition 5.2 in \cite{cabre2005parameterization}, we have that we can find a non-trivial periodic solution $\Delta$ of 
\[ (\Lo^* + \eta )\Delta = 0,\]
if and only if $e^{\eta}$ is an eigenvalue of the monodromy matrix $\Psi_T$ of system $\Lo^* \Delta =0$. Conversely, if $e^{\eta}$ is not an eigenvalue of the monodromy matrix $\Psi_T$, then given any periodic function $R$, there exists a unique periodic function $\Delta$ solving the equation
\[(\Lo^*+\eta) \Delta =R. \]

Thus, adjoint equations \eqref{eq:NTorder0}-\eqref{eq:NSorder0} amount to $\nabla \Theta^{(0)}:=\nabla \Theta (\gamma(\theta))$, $\nabla \Sigma_j^{(0)}=\nabla \Sigma_j (\gamma(\theta))$ being periodic eigenfunctions of $\Lo^*$ with eigenvalue $\lambda_j$, $j=0,\ldots, d-1$ (recall $\lambda_0=0$). So, the equations \eqref{eq:NTorder0}-\eqref{eq:NSorder0} can be written as
\[(\Lo^* - \lambda_j) \Delta = 0. \]

Therefore, the equations can then be solved if $e^{-\lambda_j}$ is an eigenvalue of the monodromy matrix $\Psi_T$ and the solution is given by expression in \eqref{eq:fsa}. Moreover, the solution is unique except for a free mutiplicative factor that is determined according to the normalization conditions \eqref{eq:ncT0} and \eqref{eq:ncS0}.

The rest of the equations in \eqref{eq:NTordern}-\eqref{eq:NSordern} are of the form
\[(\Lo^*+m\lambda_s) \Delta = R, \quad m \geq 0,\]
which can be solved uniquely if $-m \lambda_s$ is not in the spectrum of $\Lo^*$, which happens if an only if $e^{m \lambda_s}$ is not an eigenvalue of the monodromy matrix $\Psi_T$. Since the eigenvalues of the monodromy matrix $\Psi_T$ are $e^{-\overline \lambda_j}$ (see Section~\ref{sec:app-relation-eigen}) and $Re(\lambda_j)<0$ except for $\lambda_0=0$, therefore we only have problems when $m=0$. This situation happens in equation \eqref{eq:NTorder1}, for which we need to apply the normalization condition \eqref{eq:Nc-Wilson-v2}. For the rest of cases ($m \geq 1$) the equation can be solved uniquely.

\subsection{Relationship between the eigenvalues of the monodromy matrices $\Phi_T$ and $\Psi_T$}\label{sec:app-relation-eigen}

In this section, we provide the relationship between the eigenvalues of the monodromy matrix $\Phi_T$ of system \eqref{eq:mjVarEqs} (which are denoted $e^{\lambda_j}$) and the eigenvalues of the monodromy matrix $\Psi_T$ of system \eqref{eq:NTorder0}. 

Recall that given the operator
\[\Lo:=\frac{1}{T} \frac{d}{d \theta} - DX(K_{\mathbf{0}} (\theta)),\]
then $-\Lo^*$ is its adjoint operator, i.e.
$(u,\Lo u)=(-\Lo^*u,v)$, where $(\cdot,\cdot)$ denotes the standard inner product on the Hilbert space of $T$-periodic functions in $\mathbb{R}^d$ (see \cite{ErmentroutTerman2010}). Therefore, we have the following relationship for the spectrum of both operators 
\[\sigma(-\Lo^*)=\overline{\sigma(\Lo)}:=\{z\in \mathbb {C} :{\overline {z}}\in \sigma (\Lo)\}.\]

Equation \eqref{eq:eqMione} writes as
\[(\Lo + \lambda_j)K_{e_j}=0,\]
therefore, $-\lambda_j$ is an eigenvalue in the spectrum of $\Lo$, so $-\overline \lambda_j$ is in the spectrum of $-\Lo^*$, i.e.
\[ (- \Lo^* + \overline \lambda_j) \Delta =0\]
for some non-trivial periodic function $\Delta$. Notice that this is equivalent to 
\[ ( \Lo^* - \overline \lambda_j) \Delta =0,\]
and we have that $\overline \lambda_j$ is an eigenvalue of $\Lo^*$, or equivalently $e^{-\overline \lambda_j}$ is an eigenvalue of the monodromy matrix $\Psi_T$.

\bibliographystyle{abbrv}
\bibliography{references}

\begin{thebibliography}{10}

\bibitem{bonnin17}
M.~Bonnin.
\newblock Amplitude and phase dynamics of noisy oscillators.
\newblock {\em International Journal of Circuit Theory and Applications},
  45(5):636--659, 2017.

\bibitem{cabre2003parameterization}
X.~Cabr{\'e}, E.~Fontich, and R.~de~la Llave.
\newblock The parameterization method for invariant manifolds {I}: manifolds
  associated to non-resonant subspaces.
\newblock {\em Indiana University mathematics journal}, pages 283--328, 2003.

\bibitem{cabre2003parameterization2}
X.~Cabr{\'e}, E.~Fontich, and R.~de~la Llave.
\newblock The parameterization method for invariant manifolds {II}: regularity
  with respect to parameters.
\newblock {\em Indiana University mathematics journal}, pages 329--360, 2003.

\bibitem{cabre2005parameterization}
X.~Cabr{\'e}, E.~Fontich, and R.~De~La~Llave.
\newblock The parameterization method for invariant manifolds {III}: overview
  and applications.
\newblock {\em Journal of Differential Equations}, 218(2):444--515, 2005.

\bibitem{castejon2013phase}
O.~Castej{\'o}n, A.~Guillamon, and G.~Huguet.
\newblock Phase-amplitude response functions for transient-state stimuli.
\newblock {\em The Journal of Mathematical Neuroscience}, 3(1):13, 2013.

\bibitem{castelli2015parameterization}
R.~Castelli, J.-P. Lessard, and J.~D. Mireles~James.
\newblock Parameterization of invariant manifolds for periodic orbits {I}:
  Efficient numerics via the floquet normal form.
\newblock {\em SIAM Journal on Applied Dynamical Systems}, 14(1):132--167,
  2015.

\bibitem{Cestnik2022}
R.~Cestnik, E.~T.~K. Mau, and M.~Rosenblum.
\newblock {Inferring oscillator's phase and amplitude response from a scalar
  signal exploiting test stimulation}.
\newblock {\em New Journal of Physics}, 24(12):123012, 2022.

\bibitem{Clusella2023}
P.~Clusella, G.~Deco, M.~L. Kringelbach, G.~Ruffini, and J.~Garcia-Ojalvo.
\newblock Complex spatiotemporal oscillations emerge from transverse
  instabilities in large-scale brain networks.
\newblock {\em PLOS Computational Biology}, 19(4):e1010781, 2023.

\bibitem{dumont2019macroscopic}
G.~Dumont and B.~Gutkin.
\newblock Macroscopic phase resetting-curves determine oscillatory coherence
  and signal transfer in inter-coupled neural circuits.
\newblock {\em PLoS computational biology}, 15(5):e1007019, 2019.

\bibitem{ErmentroutTerman2010}
B.~Ermentrout and D.~Terman.
\newblock {\em {Mathematical foundations of neuroscience.}}
\newblock New York : Springer, 2010.

\bibitem{ErmentroutKopell91}
G.~B. Ermentrout and N.~Kopell.
\newblock Multiple pulse interactions and averaging in systems of coupled
  neural oscillators.
\newblock {\em J. Math. Biol.}, 29(3):195--217, 1991.

\bibitem{griewank2008evaluating}
A.~Griewank and A.~Walther.
\newblock {\em Evaluating derivatives: principles and techniques of algorithmic
  differentiation}, volume 105.
\newblock Siam, 2008.

\bibitem{guckenheimer1975}
J.~Guckenheimer.
\newblock Isochrons and phaseless sets.
\newblock {\em Journal of Mathematical Biology}, 1(3):259--273, 1975.

\bibitem{guckenheimer2013nonlinear}
J.~Guckenheimer and P.~Holmes.
\newblock {\em Nonlinear oscillations, dynamical systems, and bifurcations of
  vector fields}, volume~42.
\newblock Springer Science \& Business Media, 2013.

\bibitem{guillamon2009computational}
A.~Guillamon and G.~Huguet.
\newblock A computational and geometric approach to phase resetting curves and
  surfaces.
\newblock {\em SIAM Journal on Applied Dynamical Systems}, 8(3):1005--1042,
  2009.

\bibitem{haro2016}
{\`A}.~Haro, M.~Canadell, J.-L. Figueras, A.~Luque, and J.-M. Mondelo.
\newblock {\em The Parameterization Method for Invariant Manifolds}.
\newblock Springer, 2016.

\bibitem{hoppensteadt2012}
F.~C. Hoppensteadt and E.~M. Izhikevich.
\newblock {\em Weakly connected neural networks}, volume 126.
\newblock Springer Science \& Business Media, 2012.

\bibitem{huguet2013computation}
G.~Huguet and R.~de~la Llave.
\newblock Computation of limit cycles and their isochrons: fast algorithms and
  their convergence.
\newblock {\em SIAM Journal on Applied Dynamical Systems}, 12(4):1763--1802,
  2013.

\bibitem{Huguet2012}
G.~Huguet, R.~de~la Llave, and Y.~Sire.
\newblock Computation of whiskered invariant tori and their associated
  manifolds: New fast algorithms.
\newblock {\em Discrete and Continuous Dynamical Systems - A},
  32(4):1309–1353, 2012.

\bibitem{kaiser2021data}
E.~Kaiser, J.~N. Kutz, and S.~L. Brunton.
\newblock Data-driven discovery of koopman eigenfunctions for control.
\newblock {\em Machine Learning: Science and Technology}, 2(3):035023, 2021.

\bibitem{Kralemann2014}
B.~Kralemann, A.~Pikovsky, and M.~Rosenblum.
\newblock {Reconstructing effective phase connectivity of oscillator networks
  from observations}.
\newblock {\em New Journal of Physics}, 16(8):085013, aug 2014.

\bibitem{Krauskopfetal05}
B.~Krauskopf, H.~M. Osinga, E.~J. Doedel, M.~E. Henderson, J.~Guckenheimer,
  A.~Vladimirsky, M.~Dellnitz, and O.~Junge.
\newblock A survey of methods for computing (un)stable manifolds of vector
  fields.
\newblock {\em International Journal of Bifurcation and Chaos},
  15(03):763--791, 2005.

\bibitem{kuramoto2003chemical}
Y.~Kuramoto.
\newblock {\em Chemical oscillations, waves, and turbulence}.
\newblock Courier Corporation, 2003.

\bibitem{mauroy2018global}
A.~Mauroy and I.~Mezi{\'c}.
\newblock Global computation of phase-amplitude reduction for limit-cycle
  dynamics.
\newblock {\em Chaos: An Interdisciplinary Journal of Nonlinear Science},
  28(7):073108, 2018.

\bibitem{mauroy2013isostables}
A.~Mauroy, I.~Mezi{\'c}, and J.~Moehlis.
\newblock Isostables, isochrons, and koopman spectrum for the action--angle
  representation of stable fixed point dynamics.
\newblock {\em Physica D: Nonlinear Phenomena}, 261:19--30, 2013.

\bibitem{Nicks2024}
R.~Nicks, R.~Allen, and S.~Coombes.
\newblock Insights into oscillator network dynamics using a phase-isostable
  framework.
\newblock {\em Chaos: An Interdisciplinary Journal of Nonlinear Science},
  34(1), 2024.

\bibitem{Orieux2024}
M.~Orieux, A.~Guillamon, and G.~Huguet.
\newblock Optimal control of oscillatory neuronal models with applications to
  communication through coherence.
\newblock {\em Physica D: Nonlinear Phenomena}, 467:134267, 2024.

\bibitem{osinga2010continuation}
H.~M. Osinga and J.~Moehlis.
\newblock Continuation-based computation of global isochrons.
\newblock {\em SIAM Journal on Applied Dynamical Systems}, 9(4):1201--1228,
  2010.

\bibitem{perezrole}
A.~P{\'e}rez~Cervera.
\newblock {\em On the role of oscillatory dynamics in neural communication}.
\newblock PhD thesis, Universitat Polit{\`e}cnica de Catalunya, 2019.

\bibitem{PIK01}
A.~Pikovsky, M.~G. Rosenblum, and J.~Kurths.
\newblock {\em Synchronization, A Universal Concept in Nonlinear Sciences}.
\newblock Cambridge University Press, 2001.

\bibitem{PerezCervera2020}
A.~Pérez-Cervera, T.~M-Seara, and G.~Huguet.
\newblock Global phase-amplitude description of oscillatory dynamics via the
  parameterization method.
\newblock {\em Chaos: An Interdisciplinary Journal of Nonlinear Science},
  30(8), 2020.

\bibitem{PerezSH20}
A.~Pérez-Cervera, T.~M. Seara, and G.~Huguet.
\newblock Phase-locked states in oscillating neural networks and their role in
  neural communication.
\newblock {\em Communications in Nonlinear Science and Numerical Simulation},
  80:104992, 2020.

\bibitem{ReynerHuguet22}
D.~Reyner-Parra and G.~Huguet.
\newblock Phase-locking patterns underlying effective communication in exact
  firing rate models of neural networks.
\newblock {\em PLOS Computational Biology}, 18(5):1--41, 05 2022.

\bibitem{Rosenblum2001}
M.~G. Rosenblum and a.~S. Pikovsky.
\newblock {Detecting direction of coupling in interacting oscillators.}
\newblock {\em Physical review. E, Statistical, nonlinear, and soft matter
  physics}, 64(4 Pt 2):045202, 2001.

\bibitem{Segneri2020}
M.~Segneri, H.~Bi, S.~Olmi, and A.~Torcini.
\newblock Theta-nested gamma oscillations in next generation neural mass
  models.
\newblock {\em Frontiers in Computational Neuroscience}, 14, 2020.

\bibitem{shirasaka2017phase}
S.~Shirasaka, W.~Kurebayashi, and H.~Nakao.
\newblock Phase-amplitude reduction of transient dynamics far from attractors
  for limit-cycling systems.
\newblock {\em Chaos: An Interdisciplinary Journal of Nonlinear Science},
  27(2):023119, 2017.

\bibitem{simo1990analytical}
C.~Sim{\'o}.
\newblock On the analytical and numerical approximation of invariant manifolds.
\newblock In {\em Les M{\'e}thodes Modernes de la M{\'e}canique C{\'e}leste.
  Modern methods in celestial mechanics}, pages 285--329, 1990.

\bibitem{strogatzbook}
S.~H. Strogatz.
\newblock {\em {Nonlinear Dynamics and Chaos: With Applications to Physics,
  Biology, Chemistry and Engineering}}.
\newblock Westview Press, 1994.

\bibitem{Taher2020}
H.~Taher, A.~Torcini, and S.~Olmi.
\newblock Exact neural mass model for synaptic-based working memory.
\newblock {\em PLOS Computational Biology}, 16(12):e1008533, 2020.

\bibitem{James16}
J.~B. van~den Berg, J.~D. Mireles~James, and C.~Reinhardt.
\newblock Computing (un)stable manifolds with validated error bounds:
  non-resonant and resonant spectra.
\newblock {\em J. Nonlinear Sci.}, 26(4):1055--1095, 2016.

\bibitem{wedgwood2013phase}
K.~C. Wedgwood, K.~K. Lin, R.~Thul, and S.~Coombes.
\newblock Phase-amplitude descriptions of neural oscillator models.
\newblock {\em The Journal of Mathematical Neuroscience}, 3(1):2, 2013.

\bibitem{williams2015data}
M.~O. Williams, I.~G. Kevrekidis, and C.~W. Rowley.
\newblock A data--driven approximation of the koopman operator: Extending
  dynamic mode decomposition.
\newblock {\em Journal of Nonlinear Science}, 25:1307--1346, 2015.

\bibitem{Wilson2019}
D.~Wilson.
\newblock Isostable reduction of oscillators with piecewise smooth dynamics and
  complex floquet multipliers.
\newblock {\em Physical Review E}, 99(2), 2019.

\bibitem{Wilson2020b}
D.~Wilson.
\newblock {A data-driven phase and isostable reduced modeling framework for
  oscillatory dynamical systems}.
\newblock {\em Chaos}, 30(1), 2020.

\bibitem{Wilson2020}
D.~Wilson.
\newblock Phase-amplitude reduction far beyond the weakly perturbed paradigm.
\newblock {\em Physical Review E}, 101(2), 2020.

\bibitem{Wilson2023}
D.~Wilson.
\newblock {A direct method approach for data-driven inference of high accuracy
  adaptive phase-isostable reduced order models}.
\newblock {\em Physica D: Nonlinear Phenomena}, 446:133675, 2023.

\bibitem{wilson2018greater}
D.~Wilson and B.~Ermentrout.
\newblock Greater accuracy and broadened applicability of phase reduction using
  isostable coordinates.
\newblock {\em Journal of mathematical biology}, 76(1-2):37--66, 2018.

\bibitem{moehliswilsonpre2016}
D.~Wilson and J.~Moehlis.
\newblock Isostable reduction of periodic orbits.
\newblock {\em Physical Review E}, 94(5):052213, 2016.

\bibitem{winfree1967biological}
A.~T. Winfree.
\newblock Biological rhythms and the behavior of populations of coupled
  oscillators.
\newblock {\em Journal of theoretical biology}, 16(1):15--42, 1967.

\bibitem{winfree2001geometry}
A.~T. Winfree.
\newblock {\em The geometry of biological time}, volume~12.
\newblock Springer Science \& Business Media, 2001.

\bibitem{Yeldesbay2019}
A.~Yeldesbay, G.~R. Fink, and S.~Daun.
\newblock {Reconstruction of effective connectivity in the case of asymmetric
  phase distributions}.
\newblock {\em Journal of Neuroscience Methods}, 317:94--107, 2019.

\bibitem{Yeldesbay2024}
A.~Yeldesbay, G.~Huguet, and S.~Daun.
\newblock Reconstruction of phase-amplitude dynamics from electrophysiological
  signals.
\newblock Preprint arXiv:2406.05073, 2024.

\end{thebibliography}

\end{document}